\documentclass{article}

\usepackage{amsmath,amsthm}
\usepackage{amssymb}
\usepackage{graphics}
\usepackage{color}
\usepackage{xcolor}
\usepackage{amsfonts}
\usepackage{graphicx}
\usepackage{epstopdf}
\usepackage{algorithmic}
\newtheorem{prop}{Proposition}

\usepackage{PRIMEarxiv}

\usepackage[utf8]{inputenc} 
\usepackage[T1]{fontenc}    
\usepackage{hyperref}       
\usepackage{url}            
\usepackage{booktabs}       
\usepackage{nicefrac}       
\usepackage{fancyhdr}       
\graphicspath{{media/}}     

\newtheorem{lemma}{Lemma}

\title{Generalized weakly corrected Milstein solutions to stochastic differential equations}

\author{Tapas Tripura\\
  Department of Applied Mechanics\\
  Indian Institute of Technology Delhi\\
  \texttt{tapas.t@am.iitd.ac.in} \\
  
  \And
  Budhaditya Hazra \\
  Department of Civil Engineering\\
  Indian Institute of Technology Guwahati\\
  \texttt{budhaditya.hazra@iitg.ac.in} \\
  
  \And
      Souvik Chakraborty \\
  Department of Applied Mechanics\\
  School of Artificial Intelligence (ScAI)\\
  Indian Institute of Technology Delhi\\
  \texttt{souvik@am.iitd.ac.in} \\
}

\begin{document}
\maketitle

\begin{abstract}
  In this work, weakly corrected explicit, semi-implicit and implicit Milstein approximations are presented for the solution of nonlinear stochastic differential equations. The solution trajectories provided by the Milstein schemes are corrected by employing the \textit{change of measures}, aimed at removing the error associated with the diffusion process incurred due to the transformation between two probability measures. The change of measures invoked in the Milstein schemes ensure that the solution from the mapping is measurable with respect to the filtration generated by the error process. The proposed scheme incorporates the error between the approximated mapping and the exact representation as an innovation, that is accounted for, in the Milstein trajectories as an additive term. Numerical demonstration using a parametrically and non-parametrically excited stochastic oscillators, demonstrates the improvement in the solution accuracy for the corrected schemes with coarser time steps when compared with the classical Milstein approximation with finer time steps.
\end{abstract}


\section{Introduction}
Stochastic differential equations (SDEs) are of significant importance due to their ability to account for the randomness in prediction of the response evolution of the dynamical system under consideration \cite{calin2015informal,klebaner2005introduction}. In  recently published works the author has demonstrated how SDEs can be utilized in a unified framework with Health Monitoring studies to constitute a new perspective for application of SDEs \cite{tripura2020real,bhowmik2019first}. A SDE is driven by its deterministic component under the effect of additive erratic component, these are referred as drift and diffusion. Both the components are allowed to depend firstly, on time and secondly, on systems states simultaneously. Towards this, considering ${{\bf{X}}}_t^k=({\rm{X}}_t^1, \ldots {\rm{X}}_t^m)$ to be a $m$-dimensional $n$-factor SDE where each of the $m$-diffusions is driven by $n$ Wiener processes (${\bf W}_j(t),j = 1, \ldots n$).  
\begin{equation}\label{sdegen}
\begin{array}{l}
d{\bf{X}}_t = {{\bf g}}\left( {t,{{\bf{X}}_t}} \right)dt + \sum\limits_{j = 1}^n {{\bf f}_j\left( {t,{{\bf{X}}_t}} \right)} d{{\bf W}_j}\left( t \right); \quad
{\bf X}(t=t_0)={\bf X}_0; \quad t \in [0,T]
\end{array}
\end{equation}
Let the above SDE is defined under the probability space $\left( {\Omega ,\mathcal{F},{{({\mathcal{F}_t})}_t},P} \right)$, with natural filtration $({{({\mathcal{F}_t})},0 \le t \le T)}$ constructed from sub $\sigma$-algebras of $\mathcal{F}$. Then for each of the $m$ diffusions, ${\bf{X}}_t \in {\mathbb{R}^m}$ denotes the ${{\mathcal{F}_t}}$-measurable state vector, ${{\bf g}}\left( {t,{{\bf{X}}_t}} \right) \in {\mathbb{R}^m}$ is drift vector, ${{\bf f}}\left( {t,{{\bf{X}}_t}} \right) \in {\mathbb{R}^{m \times n}}$ is a volatility coefficient matrix and ${{\bf W}_j}\left( t \right) \in {\mathbb{R}^n}$ is the Wiener process with respect to the probability measure P. For $m$=$n$=1, the one dimensional one factor SDE can be written as:
\begin{equation}\label{sdeg}
\begin{array}{l}
d{{\bf{X}}_t} = {\bf g}\left( {t,{{\bf {X}}_t}} \right)dt + {\bf f}\left( {t,{{\bf {X}}_t}} \right)d{\bf W}\left( t \right); \quad {\bf X}(t=t_0)={\bf X}_0; \quad t \in [0,T]
\end{array}
\end{equation}
Then there exists a unique solution for the above SDE if \cite{hassler2016stochastic}: \textit{firstly}, ${\bf g}\left( {t,{{\bf {X}}_t}} \right)$ and ${\bf f}\left( {t,{{\bf {X}}_t}} \right)$ must be partially differentiable with respect to $X(t)$ and also continuous in $X(t)$, in other words it demands that there exists a constant ${\mathbb{C}}$ such that they are Lipschitz continuous: $\left| {{\bf{g}}\left( {t,{{\bf{X}}_t}} \right)} - {{\bf{g}}\left( {t,{{\bf{Y}}_t}} \right)} \right| + \left| {{\bf{f}}\left( {t,{{\bf{X}}_t}} \right)} - {{\bf{f}}\left( {t,{{\bf{Y}}_t}} \right)} \right| \le {\mathbb{C}}\left| {\bf{X}} - {\bf{Y}} \right|$, \textit{secondly}, the growth of the diffusion process is restricted or there exists two ${\mathbb{C}_1}$ and ${\mathbb{C}_2}$ such that $\left| {{\bf{g}}\left( {t,{{\bf{X}}_t}} \right)} \right| + \left| {{\bf{f}}\left( {t,{{\bf{X}}_t}} \right)} \right| \le {\mathbb{C}_1} + {\mathbb{C}_2}\left| {\bf{X}} \right|$, \textit{thirdly}, it requires a well defined initial point ${{\bf{X}}_{{t_0}}}$ independent of ${\bf{W}}(t)$ and $E\left( {{\bf{X}}_{{t_0}}^2} \right) < \infty $. In absence of analytical solutions the solution of SDEs are attempted using numerical schemes of desired accuracy. One of the popular approaches is the Taylor expansion of Ito integrals \cite{kloeden1992higher}. The Ito-Taylor expansion of the Eq. (\ref{sdeg}) is derived by substituting the SDE into Ito's-lemma for differentials and then expanding using the Taylor series \cite{tripura2020ito}. For a sufficiently smooth function ${b}\left( {{t},{{\bf{X}}_{t}}} \right)$, the Ito-Taylor expansion between $t \in [{t_0},t + h]$ is expressed as,
\begin{equation}\label{Milstein}
\begin{array}{ll}
b\left( {{t_{i + 1}},{{\bf{X}}_{{t_{i + 1}}}}} \right) &= b\left( {{t_0},{{\bf{X}}_{{t_0}}}} \right) + {\Im ^0}b\left( {t,{{\bf{X}}_t}} \right)\int\limits_t^{{t_{i + 1}}} {ds}  + {\Im ^1}b\left( {t,{{\bf{X}}_t}} \right)\int\limits_t^{{t_{i + 1}}} {d{W_s}}  + \\& {\Im ^1}{\Im ^1}b\left( {t,{{\bf{X}}_t}} \right)\int\limits_t^{{t_{i + 1}}} {\int\limits_t^{{s_1}} {d{W_{{s_2}}}} d{W_{{s_1}}}}  + \ldots +\\& {\Im ^{{j_1}}}{\Im ^{{j_2}}} \ldots {\Im ^{{j_k}}}b\left( {t,{{\bf{X}}_t}} \right)\int\limits_t^{{t_{i + 1}}} {\int\limits_t^{{s_1}} { \ldots \int\limits_t^{{s_{k - 1}}} {dP_{{s_k}}^{{j_k}} \ldots } dP_{{s_2}}^{{j_k}}dP_{{s_1}}^{{j_k}}} }  +  \ldots  + \\&
\int\limits_t^{t + h} {\int\limits_t^{{s_1}} { \cdots \int\limits_t^{{s_{k - 1}}} {{\Im ^{{j_k}}}{\Im ^{{j_{k - 1}}}} \ldots {\Im ^{{j_1}}}b\left( {{s_k},{\bf{X}}\left( {{s_k}} \right)} \right)dP_{{s_k}}^{{j_k}}}  \cdots dP_{{s_2}}^{{j_k}}dP_{{s_1}}^{{j_k}}} } 
\end{array}
\end{equation}
where ${\Im ^0}(.)$ and ${\Im ^1}(.)$ are the stochastic moments \cite{tripura2020ito}, this moments helps in generating the corresponding discrete form of an Ito-Taylor numerical scheme. In the truncated Ito-Taylor direct integration schemes, two types of integrals are observed, these are ${I_{{j_1}{j_2} \ldots {j_k}}\left( \Delta t \right)}$ and $ {I_{{j_1}{j_2} \ldots {j_k}}\left( {g,\Delta t} \right)} $, the first one appear in the corresponding truncated Ito-Taylor numerical schemes, while the other one arises in the remainders of the Ito-Taylor expansion. In the generic form these stochastic integrals are \cite{milstein1998balanced},
\begin{equation}
\label{type}
\resizebox{\textwidth}{!}{$
	\begin{array}{l}
	{I_{{j_1}{j_2} \ldots {j_k}}}\left( {\Delta t} \right) = \int\limits_t^{t + h} {\int\limits_t^{{s_1}} { \cdots \int\limits_t^{{s_{k - 1}}} {d{P_{{j_k}}}\left( {{s_k}} \right)}  \cdots d{P_{{j_2}}}\left( {{s_2}} \right)} } d{P_{{j_1}}}\left( {{s_1}} \right) \\
	{I_{{j_1}{j_2} \ldots {j_k}}}\left( {g,\Delta t} \right) = \int\limits_t^{t + h} {\int\limits_t^{{s_1}} { \cdots \int\limits_t^{{s_{k - 1}}} {{\Im ^{{j_k}}}{\Im ^{{j_{k - 1}}}} \ldots {\Im ^{{j_1}}}g\left( {{s_k},X_{{s_k}}} \right)d{P_{{j_k}}}\left( {{s_k}} \right)}  \cdots d{P_{{j_2}}}\left( {{s_2}} \right)} } d{P_{{j_1}}}\left( {{s_1}} \right)\\
	\end{array}
	$}
\end{equation}
For ${j_k}$ = 0, $d{P_{{j_k}}}$ becomes $ds$, and for ${j_k} \ne$ 0, $d{P_{{j_k}}}$ can be substituted by $dW$ \cite{milstein1998balanced}. There are various numerical techniques such as Euler-Maruyama \cite{maruyama1955continuous}, Milstein \cite{milstein1998balanced}, stochastic versions of Heun \cite{burrage2007numerical}, Runge-Kutta \cite{ruemelin1982numerical} and Newmark methods, Strong Taylor 1.5 \cite{kloeden1992higher}, Weak Taylor of order $O(\Delta t^{3})$ \cite{tripura2020ito}, which are utilized in Monte Carlo (MC) framework, further can be categorized under strong and weak schemes. Both strong and weak schemes exist either as explicit, or as implicit \cite{milstein1998balanced}, or a combination of both \cite{omar2011composite}. The implicit schemes are more robust when numerical stability is consider than explicit one, however, is limited when coarser time steps are considered. These schemes are formulated by truncating terms of required degree from the Ito-Taylor expansion \cite{kloeden1992higher,milstein1998balanced} which are either mathematically involved (i.e. higher order schemes) or  computationally exhaustive (i.e. lower order schemes). The strong schemes in this regard are characterized by their ability to estimate pathwise dynamical response, while the weak schemes tend to obtain the statistical moments of functions of dynamical response using random variables of simpler distributions rather then MSIs \cite{oksendal2013stochastic}. The application of the strong schemes are limited to small degree-of-freedom system due to the involvement in mathematical derivations, and the weak ones are although computationally efficient with identical $\Delta t$ but are not robust for pathwise system response prediction thereby not implementable to control and estimation problems. The interest of the present work revolves around the formulation of a weak corrections based on the change of measures (Girsanov transformation) \cite{girsanov1960transforming,oksendal2013stochastic,ogawa2009importance} for the explicit, semi-implicit and implicit versions of the Milstein scheme \cite{omar2011composite} towards counter of the limitations invoked by the mathematical involvement of strong and higher order numerical schemes and at the same time facilitating relaxation in requirement of finer time steps.

Two probability measure is utilized in the proposed weakly corrected schemes where a correction factor is obtained from the filtration arising in the transformation between the measures. The correction in this study is additive in nature \cite{sarkar2016weakly} and idealized as a Radon-Nikodym derivative, which is also a solution to a scalar SDE in exponential form \cite{kanjilal2017girsanov,dereudre2017exact}. It defines the distribution of weights among the ensembles of the linearized system responses and performs the sampling of linearized solution based on the weights \cite{sarkar2016weakly,raveendran2013nearly}. The stochastic exponentials arising due to the MSIs involved in exponential Radon-Nikodym derivative tend to hamper the estimation of the proper weight/fitness of the ensembles \cite{raveendran2013nearly,baldi2017introduction}, thus the correction is incorporated into the Milstein solution as an additive term in the current framework \cite{sarkar2016weakly}. To derive an appropriate additive correction term for the Milstein schemes, the central idea is the concept of change of measure, such that the Milstein approximated integrated process is measurable with respect to the filtration generated by the error, that will in turn improve the convergence of the Milstein schemes \cite{kushner1967dynamical,kallianpur1969stochastic}. The innovation thus is obtained as an integration error and is used for updating the predicted realizations of the solution using present observations. The fidelity of the proposed schemes are judged in the light of three nonlinear oscillators in the present study.

Rest of The paper is arranged as follows: \textbf{Section 2}: a short background on types of Milstein schemes are provided. \textbf{Section 3}: the proposed Girsanov corrected Milstein schemes are briefly discussed with sufficient mathematical deductions. \textbf{Section 4}: numerical illustration using a fairly representative class of non-linear mechanical oscillators excited using Gaussian white noise are provided. \textbf{Section 5}: the paper is concluded by highlighting key achievements of the work.

\section{Milstein schemes and Radon-Nikodym derivative for Wiener process}

\subsection{Milstein schemes}\label{milstein_schemes}
The Milstein approximation to the time evolution of the SDE (\ref{sdeg}) is obtained by truncating the Ito-Taylor expansion (\ref{Milstein}) at the Wiener integral of multiplicity 1.0. The three Milstein approximations over the interval, $\Delta t = ({t_i} - {t_{i - 1}})$, in the Ito sense are given as \cite{omar2011composite}:
\begin{enumerate}
	\item Explicit Milstein scheme (ML):
	\begin{equation}\label{explicit}
	\begin{array}{ll}
	{{\bf{X}}_{{t_i}}} =& {{\bf{X}}_{{t_{i - 1}}}} + {\bf g}\left( {{t_{i - 1}},{{\bf X}_{{t_{i - 1}}}}} \right)\Delta t + {\bf f}\left( {{t_{i - 1}},{{\bf X}_{{t_{i - 1}}}}} \right)\Delta {{\bf W}_{{t_i}}}  
	+ \\& \qquad
	\frac{1}{2}{\bf f}\left( {{t_{i - 1}},{{\bf X}_{{t_{i - 1}}}}} \right){\mathbf{f}}{\left( {{t_{i - 1}},{{\mathbf{X}}_{{t_{i - 1}}}}} \right)^T}\left( {\Delta {{\bf W}_{{t_i}}}^2 - \Delta t} \right)
	\end{array}
	\end{equation}
	\item Semi-implicit Milstein scheme (SIML):
	\begin{equation}\label{semi_implicite_sde}
	\begin{array}{ll}
	{{\bf{X}}_{{t_i}}} =& {{\bf{X}}_{{t_{i - 1}}}} + {\bf g}\left( {{t_{i}},{{\bf X}_{{t_{i}}}}} \right)\Delta t + {\bf f}\left( {{t_{i - 1}},{{\bf X}_{{t_{i - 1}}}}} \right)\Delta {{\bf W}_{{t_i}}} + \\& \qquad
	\frac{1}{2}{\bf f}\left( {{t_{i - 1}},{{\bf X}_{{t_{i - 1}}}}} \right){\mathbf{f}}{\left( {{t_{i - 1}},{{\mathbf{X}}_{{t_{i - 1}}}}} \right)^T}\left( {\Delta {{\bf W}_{{t_i}}}^2 - \Delta t} \right)
	\end{array}
	\end{equation}
	\item Implicit Milstein Scheme (IML):
	\begin{equation}\label{implicit}
	\begin{array}{ll}
	{{\bf{X}}_{{t_i}}} =& {{\bf{X}}_{{t_{i - 1}}}} + {\bf g}\left( {{t_{i}},{{\bf X}_{{t_{i}}}}} \right)\Delta t + {\bf f}\left( {{t_{i}},{{\bf X}_{{t_{i}}}}} \right)\Delta {{\bf W}_{{t_i}}} + \\& \qquad
	\frac{1}{2}{\bf f}\left( {{t_{i}},{{\bf X}_{{t_{i}}}}} \right){\mathbf{f}}{\left( {{t_{i - 1}},{{\mathbf{X}}_{{t_{i - 1}}}}} \right)^T}\left( {\Delta {{\bf W}_{{t_i}}}^2 - \Delta t} \right)
	\end{array}
	\end{equation}
\end{enumerate}
where, ${}^m{{\bf{X}}_{{t}}} \in {\mathbb{R}^m}$ represents the corresponding Miltein based approximated solution within the time increment $\Delta t$=$(t_i - t_{i-1})$ and Wiener increment $\Delta {\bf W}$=$({\bf W}_{t_i}-{\bf W}_{t_{i-1}})$. In the absence of multiple Wiener integrals of higher order of smallness, the Milstein scheme attains a strong order of convergence of $O ({\Delta t}^{1.0})$.

This paper focuses on constructing the governing evolution equation of the weakly corrected form of the ${}^m{{\bf{X}}_{{t_i}}}$, mentioned above. The evolutionary equation is in the form of a conditional expectation of ${}^m{{\bf{X}}_{{t_i}}}$, that follows: ${\pi _t}({\bf{X}}) = {{\rm E}_P}\left[ {\left. {^m{{\bf{X}}_t}} \right|\mathcal{F}_t} \right]$, $\mathcal{F}_t$ being the natural filtration generated by some $\mathbb{R}^m$-valued process ${{\bf \gamma} \left( {t,{{\bf{X}}_t}} \right)}$. In the later sections, it will be seen that during the transformations within two probability measures ${\pi _t}({\bf{X}})$ primarily acts as a normalizing function. Towards obtaining the governing evolutionary form, the concepts of non-linear filtering theory is being utilized in the present work. Thus to be in consistent with the nonlinear filtering theory a scalar valued function $\Phi ({\bf{X}})$ is defined for the process ${}^m{{\bf{X}}_{t}}$. Thus, the final filtering equation for the governing weakly corrected solutions will obtained as $\pi (\Phi ({\bf{X}}))$, that is remain to be shown.

\subsection{Change of measure and Radon-Nikodym derivative}
Let ${\bf X}_t$ be an Ito process that satisfies the SDE (Eq. \ref{sdeg}) on the probability measure P with the complete probability space $\left( {\Omega ,\mathcal{F},{{({\mathcal{F}_t})}_t},P} \right)$. Let Q be the another probability measure on $\left( {\Omega ,{{({\mathcal{F}_t})}_t}} \right)$. The representation of Eq. \ref{sdeg} on Q-measure satisfies the following SDE, with ${\mathbf{h}}\left( {t,{{\mathbf{X}}_t}} \right) \in {\mathbb{R}^m}$ as a modified drift follows:
\begin{equation}\label{qbrown}
	d{{\mathbf{X}}_t} = {\mathbf{h}}\left( {t,{{\mathbf{X}}_t}} \right)dt + {\mathbf{f}}\left( {t,{{\mathbf{X}}_t}} \right)d{\mathbf{\tilde W}}\left( t \right)
\end{equation} 
The probability measure Q has a density ${{\bf Z}_t}$ with respect to P, often called as Radon-Nikodym derivative. The change of measure P $\to$ Q can be effected by utilizing the derivative ${{\bf Z}_t}$ as:
\begin{equation}
dQ = {{\bf Z} _t}dP
\end{equation} 
Let ${{\bf{\gamma }}\left( {s,{{\bf{X}}_s}} \right)} \in \mathbb{R}^n$ be a progressively measurable n-dimensional process on $\left( {\Omega ,{\mathcal{F}_t}} \right)$ such that the following relation holds,
\begin{equation}\label{err_gen}
{\mathbf{f}}\left( {t,{{\mathbf{X}}_t}} \right){\mathbf{\gamma }}\left( {t,{{\mathbf{X}}_t}} \right) = {\mathbf{g}}\left( {t,{{\mathbf{X}}_t}} \right) - {\mathbf{h}}\left( {t,{{\mathbf{X}}_t}} \right)
\end{equation}
Then the Radon-Nikodym derivative ${\bf Z}_t$ satisfies to be a martingale in the form, defined as (\cite{oksendal2013stochastic,baldi2017introduction}):
\begin{equation}\label{radon}
{{\bf Z}_t} = \exp \left( {\int_{{t_{i - 1}}}^{{t_i}} {{\bf \gamma} \left( {s,{{\bf{X}}_s}} \right)} d{W_s} - \cfrac{1}{2}{\int_{{t_{i - 1}}}^{{t_i}} {{\left| {{\bf \gamma}\left( {s,{{\bf{X}}_s}} \right)} \right|} }^2}ds} \right)
\end{equation}
Considering, ${\bf{W}_t} = \left({\Omega,\mathcal{F},{{({\mathcal{F}_t})}_t},{{({\bf W}_t)}_t},P} \right)$ be an n-dimensional Wiener process with respect ot the natural augmented filtration $\mathcal{F}_t^n$ under the probability measure P. An equivalent n-measurable ${\mathcal{F}_t}$-Wiener process under $Q$ is then defined as an additional diffusion process \cite{dereudre2017exact},
\begin{equation}\label{qbrownian}
{{\tilde {\bf W}}_t} = {{\bf W}_t} + \int_{{t_{i - 1}}}^{{t_i}} {\gamma \left( {s,{{\bf{X}}_s}} \right)ds}; \quad {\tilde{\bf W}}(t=t_0)=0; \quad t\ge0
\end{equation}
Here, ${{\tilde {\bf W}}_t}$ is the n-dimensional Wiener process defined on $\left( {\Omega ,{\mathcal{F}_t}} \right)$ between [0,T] under the new probability measure-Q. With the estimate of Q-Wiener increments as: $d{{{\mathbf{\tilde W}}}_t} = d{{\mathbf{W}}_t} + \gamma \left( {t,{{\mathbf{X}}_t}} \right)dt$, one can verify that Eq. \ref{qbrown} is just a replica of Eq. \ref{sdeg} with additional drift term $\gamma \left( {t,{{\mathbf{X}}_t}} \right)dt$.

\section{Weak correction for Milstein schemes}
This section provides mathematical backbone of the three improved Milstein schemes. In order to follow the brief details of the implementation of the change of measure and nonlinear filtering theory some essential intermediate concepts and/or definitions for the appropriate treatment of the MSIs are required. \textit{First}, one will require to establish the separability condition among the terms involved in the coefficient associated with MSI of strong order 1.0 in the Milstein scheme. In the regard of an one dimensional SDE this requirement is trivial, however, in multi-degree-freedom systems perhaps one should consider a K-variate SDE primarily arising in mechanical and structural systems, the separability remains to be proved. This criteria will be used later in defining an appropriate error process. \textit{Second}, few concepts related to the cubic variation among time ($\Delta t$) and Wiener increments ($\Delta W$) needs to be established, that is essential to appropriately treat the independency between the MSIs in Milstein schemes. In the present, there exist quadratic variation among the increments, however, to formulate the governing non-linear filtration theory for the proposed improved Milstein schemes these identities are of prime significance. These criteria are one of the key entitlements of the present work, which are currently not been utilized in other inline Girsanov corrected schemes. One such scheme is the Giranov corrected Euler-Maruyama. Since the in Euler-Maruyama approximation the highest strong order of MSIs is 0.5, further the MSIs are increments itself, the independency can be established using quadratic covariation only. 

The first concept is formalized as \textbf{Proposition 1}, whereas, the second one is proved as \textbf{Lemma 1}. These are as follows:

\begin{prop}[Separability of $\sum_{k = 1}^n {{{\bf{f}}^{k,j}}\left( {t,{{\bf{X}}_t}} \right)} \frac{{\partial {{\bf{f}}^{k,j}}\left( {t,{{\bf{X}}_t}} \right)}}{{\partial {{\rm{X}}_k}}}$]\label{thm1}
	
	If a dynamical system is expressible in terms of the first-order dynamics under the probability measure ($\Omega$,$F$,$P$), in the present case possibly a SDE, through the use of an adequate statespace formulation, then it holds that,
	\begin{equation}
	\sum\limits_{k = 1}^n {{{\bf{f}}^{k,j}}\left( {t,{{\bf{X}}_t}} \right)} \frac{{\partial {{\bf{f}}^{k,j}}\left( {t,{{\bf{X}}_t}} \right)}}{{\partial {{\rm{X}}_k}}} = \sum\limits_{k = 1}^n {{{\bf{f}}^{k,j}}\left( {t,{{\bf{X}}_t}} \right)} \sum\limits_{k = 1}^n {\frac{{\partial {{\bf{f}}^{k,j}}\left( {t,{{\bf{X}}_t}} \right)}}{{\partial {{\rm{X}}_k}}}} 
	\end{equation}
	where, the term in left hand side is associated with the MSI of order of smallness 1.0. It also requires that, ${{\bf{f}}\left( {t,{{\bf{X}}_t}} \right)}$ is partially differentiable and continuous in argument $X$, in its existence in Milstein approximation it also generalized as $\textbf{f}{\left( {t,{{\bf{X}}_{t}}} \right)} \in \mathbb{R}^{m \times n}$ diffusion matrix.
\end{prop}
\begin{proof}[Proposition \ref{thm1}]
	A MDOF dynamical system is considered here whose governing dynamics is expressed by the following equation:
	\begin{equation}\label{mdof}
	{\bf{M\ddot X}}(t) + {\bf{C}}\left( {{\bf{X}},{\bf{\dot X}}} \right){\bf{\dot X}}(t) + {\bf{K}}\left( {{\bf{X}},{\bf{\dot X}}} \right){\bf{X}}(t) = {\bf{f}}\left( {t,{\bf{X}},{\bf{\dot X}},{{\bf{W}}_t}} \right){\bf{\dot W}}\left( t \right)
	\end{equation}
	The states of the system is given as, ${\bf{X}} = {({X_j};j = 1,2, \ldots n)^T}$. \textbf{M} is the constant mass matrix, ${\bf{C}}\left( {{\bf{X}},{\bf{\dot X}}} \right)$ and ${\bf{K}}\left( {{\bf{X}},{\bf{\dot X}}} \right)$ are the damping and stiffness matrices, respectively, perhaps state driven for time-varying systems. ${\bf{f}}\left( {t,{\bf{X}},{\bf{\dot X}},{{\bf{W}}_t}} \right)$ is the $n \times n$-dimensional matrix denoting either the parametric excitation or intensity of stochastic force. ${\bf{\dot W}}\left( t \right) = {({{\dot W}_j}(t);j = 1,2, \ldots n)^T}$ is the $n$-dimensional zero mean Gaussian white noise and ${\bf{W}}\left( t \right)$ is the independent and identically distributed Wiener process. To obtain the first-order SDEs two new transformations are assumed for the system states: ${\bf{X}} = {({Y_{1j}};j = 1,2, \ldots m)^T}$ and ${\bf{\dot X}} = {({Y_{2j}};j = 1,2, \ldots m)^T}$. Thus, the Eq. (\ref{mdof}) takes the form,
	\begin{equation}
	\begin{array}{ll}
	{\bf{M}}{{{\bf{\dot Y}}}_2} + {\bf{C}}\left( {{{\bf{Y}}_1},{{\bf{Y}}_2}} \right){{\bf{Y}}_2} + {\bf{K}}\left( {{{\bf{Y}}_1},{{\bf{Y}}_2}} \right){{\bf{Y}}_1} = \sum\limits_{k = 1}^n {{{\bf{f}}_k}\left( {t,{{\bf{Y}}_1},{{\bf{Y}}_2},{{\bf{W}}_t}} \right){{{\bf{\dot W}}}_k}\left( t \right)}
	\end{array}
	\end{equation}
	where, ${{{\bf{f}}_k}\left( {t,{{\bf{Y}}_1},{{\bf{Y}}_2},{{\bf{W}}_t}} \right)}$ is the $k^{th}$ column of the $m \times n$ diffusion matrix and $m = 2n$. The $m$-dimensional statespace then can be identified as,
	\begin{equation}
	\begin{array}{l}
	d{Y_{1j}}(t) = {g_{1j}}(t,{{\bf{Y}}_1},{{\bf{Y}}_2})dt\\
	d{Y_{2j}}(t) = {g_{2j}}(t,{{\bf{Y}}_1},{{\bf{Y}}_2})dt + \sum\limits_{k = 1}^n {{f_{jk}}\left( {t,{{\bf{Y}}_1},{{\bf{Y}}_2}} \right){{{\bf{\dot W}}}_k}\left( t \right)} 
	\end{array}
	\end{equation}
	Here, ${Y_{1j}}$ and ${Y_{2j}}$ denotes the displacement and velocity state of the $j^{th}$ degree-of-freedom (DOF). The complete structure of the diffusion matrix can be shown as,
	\begin{equation}
	{\bf{f}} = {\left[ {\begin{array}{*{20}{c}}
			{\bf{0}}\\
			{\sum\nolimits_k^n {{{\tilde f}_{2,k}}} }\\
			{\bf{0}}\\
			{\sum\nolimits_k^n {{{\tilde f}_{4,k}}} }\\
			\vdots \\
			{\bf{0}}\\
			{\sum\nolimits_k^n {{{\tilde f}_{m,k}}} }
			\end{array}} \right]_{m \times n}} = \left[ {\begin{array}{*{20}{c}}
		0&0& \cdots &0\\
		{{f_{2,1}}}&{{f_{2,2}}}& \cdots &{{f_{2,n}}}\\
		0&0& \cdots &0\\
		{{f_{4,1}}}&{{f_{4,2}}}& \cdots &{{f_{4,n}}}\\
		\vdots & \vdots & \ddots & \vdots \\
		0&0& \cdots &0\\
		{{f_{m,1}}}&{{f_{m,2}}}& \cdots &{{f_{m,n}}}
		\end{array}} \right]
	\end{equation}
	For the proof, the diffusion for $j^{th}$ DOF is,
	\begin{equation}
	{{\bf{f}}_j}(t,{{\bf{Y}}_1},{{\bf{Y}}_2}) = \left[ {\begin{array}{*{20}{c}}
		{{f_{1j}}(t,{{\bf{Y}}_1},{{\bf{Y}}_2})}\\
		{{f_{2j}}(t,{{\bf{Y}}_1},{{\bf{Y}}_2})}
		\end{array}} \right]
	\end{equation}
	where, diffusion terms for $j$-DOF is identified as,
	\begin{equation*}
	\begin{array}{l}
	{f_{1j}}(t,{{\bf{Y}}_1},{{\bf{Y}}_2}) = {{\bf{0}}_{1 \times n}}\\
	{f_{2j}}(t,{{\bf{Y}}_1},{{\bf{Y}}_2}) = \sum\limits_{k = 1}^n {{f_{jk}}\left( {t,{{\bf{Y}}_1},{{\bf{Y}}_2}} \right)}  = {[\begin{array}{*{20}{c}}
		{{{\tilde f}_{j1}}}&{{{\tilde f}_{j2}}}& \ldots &{{{\tilde f}_{jn}}}
		\end{array}]_{1 \times n}}
	\end{array}
	\end{equation*}
	It needs to understand here that, the coefficient of the order of smallness $O(h)^{1.0}$ in the Milstein aprroximation of a SDE, has the definition \cite{milstein1998balanced},
	\begin{equation}
	{\Im ^j}(.) = \sum\limits_{k = 1}^m {{f_{kj}}(t,{{\bf{X}}_t})} \frac{\partial }{{\partial {X_k}}}(.);j = 1,2, \ldots n
	\end{equation}
	Finally, carefully observing the diffusion matrix and identifying that, \[{\bf{f}} = \left[ {\begin{array}{*{20}{c}}
		{\sum\nolimits_k^n {{{\tilde f}_{1,k}}} }&{\sum\nolimits_k^n {{{\tilde f}_{2,k}}} }&{\sum\nolimits_k^n {{{\tilde f}_{3,k}}} }& \ldots &{\sum\nolimits_k^n {{{\tilde f}_{j,k}}}  \ldots }&{\sum\nolimits_k^n {{{\tilde f}_{m,k}}} }
		\end{array}} \right];k = 1,2, \ldots n\]
	and ${\Im ^j}\left( {{{\bf{f}}_1}(t,{{\bf{X}}_t})} \right) = \sum_{k = 1}^n {{{\bf{f}}^{k,j}}\left( {t,{{\bf{X}}_t}} \right)} \frac{{\partial {{\bf{f}}^{k,j}}\left( {t,{{\bf{X}}_t}} \right)}}{{\partial {{\rm{X}}_k}}}$, it is straightforward to note that,
	\begin{equation}
	\sum\limits_{k = 1}^n {{{\bf{f}}^{k,j}}\left( {t,{{\bf{X}}_t}} \right)} \frac{{\partial {{\bf{f}}^{k,j}}\left( {t,{{\bf{X}}_t}} \right)}}{{\partial {{\rm{X}}_k}}} = \left\{ \begin{array}{ll}
	  \sum\limits_{k = 1}^n {{{\bf{f}}^{k,j}}\left( {t,{{\bf{X}}_t}} \right)} \sum\limits_{k = 1}^n {\frac{{\partial {{\bf{f}}^{k,j}}\left( {t,{{\bf{X}}_t}} \right)}}{{\partial {{\rm{X}}_k}}}};   & \text{$j$ = even} \\
	   0;  & \text{$j$ = odd}
	\end{array} \right.
	\end{equation}
	Without imposing further more conditions, the above result for a first order uni-variate 1-factor SDE follows: ${{\bf f}\left( {{t_{i - 1}},{{\bf{X}}_{{t_{i - 1}}}}} \right)}{\bf f}^\prime{\left( {{t_{i - 1}},{{\bf{X}}_{{t_{i - 1}}}}} \right)}$=${{\bf f}\left( {{t_{i - 1}},{{\bf{X}}_{{t_{i - 1}}}}} \right)} {\bf f}^\prime{\left( {{t_{i - 1}},{{\bf{X}}_{{t_{i - 1}}}}} \right)} $. $\blacksquare$
\end{proof}

\begin{lemma}[Cubic Identities]\label{pf1}
	Let $n$ be the length of a partition $P_n$ between $[0,t]$, such that ${P_n}([0,t]):{s_0} = 0 < {s_1} < {s_2} <  \ldots  < {s_i} <  \ldots {s_n} = t$. Then for $n \to \infty$, it holds that,
	\begin{enumerate}
		\item [(i)] $\sum\limits_{i = 1}^n {{{\left( {W({s_i}) - W({s_{i - 1}})} \right)}^2}\left( {{s_i} - {s_{i - 1}}} \right)} \xrightarrow{2} \int\limits_0^t {{{\left( {dW(s)} \right)}^2}dt = 0}$
		\item [(ii)] $\sum\limits_{i = 1}^n {{{\left( {W({s_i}) - W({s_{i - 1}})} \right)}^3}} \xrightarrow{2} \int\limits_0^t {{{\left( {dW(s)} \right)}^3} = 0} $
	\end{enumerate}
\end{lemma}
\begin{proof}[Proof of Lemma \ref{pf1}]
	(i) Consider the covariation as:
	\begin{equation}
	C{V_n} = \sum\limits_{i = 1}^n {{{\left( {W({s_i}) - W({s_{i - 1}})} \right)}^2}\left( {({s_i}) - ({s_{i - 1}})} \right)} 
	\end{equation}
	The proof requires the establishment of results: ${\rm E}[C{V_n}] = 0$, ${\rm E}\left[ {{{\left( {C{V_n}} \right)}^2}} \right] = 0$, and ${\rm E}\left[ {{{\left( {C{V_n}} \right)}^3}} \right] = 0$. To prove ${\rm E}[C{V_n}] = 0$ take expectation on $C{V_n}$ and apply Fubini's theorem. Further, noting $\left( {W({s_i}) - W({s_{i - 1}})} \right)~N(0,({s_i}) - ({s_{i - 1}}))$ one gets,
	\begin{equation*}
	\begin{array}{ll}
	{\rm E}[C{V_n}] &= \sum\limits_{i = 1}^n {{\rm E}\left[ {{{\left( {W({s_i}) - W({s_{i - 1}})} \right)}^2}} \right]\left( {({s_i}) - ({s_{i - 1}})} \right)} \\
	&= \sum\limits_{i = 1}^n {{{\left( {({s_i}) - ({s_{i - 1}})} \right)}^2}} \\
	&\le \underbrace {\max }_{1 \le i \le n}\left( {({s_i}) - ({s_{i - 1}})} \right)\sum\limits_{i = 1}^n {\left( {({s_i}) - ({s_{i - 1}})} \right)} \\
	&= \underbrace {\max }_{1 \le i \le n}\left( {({s_i}) - ({s_{i - 1}})} \right){V_n}(id,t)\\
	&\to 0
	\end{array}
	\end{equation*}
	where, ${V_n}(id,t) \to t$ is the absolute variation and $id$ is the identity function such that $id(t)=t$ \cite{hassler2016stochastic}. The intuition is that ${Q_n}(id,t)=\sum_{i = 1}^n {{{\left( {({s_i}) - ({s_{i - 1}})} \right)}^2}}$ contains $n$ partitions and the magnitude of $\left( {({s_i}) - ({s_{i - 1}})} \right)$ is $1/{n^2}$, thus the sum converges $\to$ 0, as $n \to \infty$. Since, ${\rm E}[C{V_n}] = 0$ it is straightforward to establish that, ${\rm MSE}(C{V_n},0) = Var(C{V_n})$, where ${\rm MSE}(X_n,X)$ is the mean squared error as norm between the sequence $X_n$ and random variable $X$ such that ${\rm MSE}(X_n,X)={\rm E}[(X_n-X)^{2}]$ \cite{calin2015informal}. Further, one may note that $X_n$ converges in mean square to $X$ as $n \to \infty$, that is formalized as ${\rm MSE}(X_n,X) \xrightarrow{2} {\rm E}[(X_n-X)^{2}]; \quad n \to \infty $.
	\begin{equation*}
	\resizebox{\textwidth}{!}{$
	\begin{array}{ll}
	Var[C{V_n}] &= \sum\limits_{i = 1}^n {Var\left[ {{{\left( {W({s_i}) - W({s_{i - 1}})} \right)}^2}} \right]{{\left( {({s_i}) - ({s_{i - 1}})} \right)}^2}} \\&
	= \sum\limits_{i = 1}^n {\left( {{\rm E}\left[ {{{\left( {W({s_i}) - W({s_{i - 1}})} \right)}^4}} \right] - {\rm E}{{\left[ {{{\left( {W({s_i}) - W({s_{i - 1}})} \right)}^2}} \right]}^2}} \right)} \left( {({s_i}) - ({s_{i - 1}})} \right) \\&
	= \sum\limits_{i = 1}^n {\left( {3{\rm E}{{\left[ {{{\left( {W({s_i}) - W({s_{i - 1}})} \right)}^2}} \right]}^2} - \left( {({s_i}) - ({s_{i - 1}})} \right)} \right)}{{\left( {({s_i}) - ({s_{i - 1}})} \right)}^2} \\&
	= 2\sum\limits_{i = 1}^n {{{\left( {({s_i}) - ({s_{i - 1}})} \right)}^4}} \\&
	\le 2\underbrace {\max }_{1 \le i \le n}\left( {({s_i}) - ({s_{i - 1}})} \right)\sum\limits_{i = 1}^n {{{\left( {({s_i}) - ({s_{i - 1}})} \right)}^3}} \\&
	= \underbrace {\max }_{1 \le i \le n}\left( {({s_i}) - ({s_{i - 1}})} \right){Q_n}(id,t){V_n}(id,t)\\&
	\to 0 
	\end{array} $}
	\end{equation*}
	The concept is similar, as $\sum_{i = 1}^n {{{\left( {({s_i}) - ({s_{i - 1}})} \right)}^4}}$ contains $n$ numbers of $\left( {({s_i}) - ({s_{i - 1}})} \right)$, each having a magnitude of $1/{n^4}$, it can be understood that $\left( {({s_i}) - ({s_{i - 1}})} \right) \to 0$ as $n \to \infty$. The third property ${\rm E}\left[ {{{\left( {C{V_n}} \right)}^3}} \right] = 0$ is implied due to zero skewness of Gaussian random variable. The proof is complete. $\blacksquare$
	
	(ii) Similar to above problem, assume,
	\begin{equation*}
	C{V_n} = \sum\limits_{i = 1}^n {{{\left( {W({s_i}) - W({s_{i - 1}})} \right)}^3}} 
	\end{equation*}
	Noting that, the skewness ${\rm E}[X^3]$ of Gaussian distribution is zero, it is obvious that,
	\begin{equation*}
	{\rm E}[C{V_n}] = \sum\limits_{i = 1}^n {{\rm E}\left[ {{{\left( {W({s_i}) - W({s_{i - 1}})} \right)}^3}} \right]}=0
	\end{equation*}
	Since, ${\rm E}[C{V_n}]=0$, it holds that, ${\rm MSE}(C{V_n},0) = Var(C{V_n})$ and it remains to prove that ${\rm MSE}(C{V_n},0) \to 0$.
	\begin{equation*}
	\begin{array}{ll}
	Var[C{V_n}] &= \sum\limits_{i = 1}^n {Var\left[ {\left( {W({s_i}) - W({s_{i - 1}})} \right)^3} \right]} \\&
	= \sum\limits_{i = 1}^n {\left( {{\rm E}\left[ {{{\left( {W({s_i}) - W({s_{i - 1}})} \right)}^6}} \right] - {\rm E}{{\left[ {{{\left( {W({s_i}) - W({s_{i - 1}})} \right)}^3}} \right]}^2}} \right)}
	\end{array}
	\end{equation*}
	At this stage, the property of Hyper-Kurtosis is to be appropriately utilized completely solve the problem. Under the result $E[C{V_n}]=0$, the Hyper-Kurtosis is defined as,
	\begin{equation*}
	{{\bf{K}}^H} = \frac{{{\rm E}\left[ {{{\left( {X - {\rm E}\left[ X \right]} \right)}^6}} \right]}}{{{{\left( {{\rm E}\left[ {{{\left( {X - {\rm E}\left[ X \right]} \right)}^2}} \right]} \right)}^3}}} \Rightarrow {\rm E}\left[ {\frac{{{{\left( X \right)}^6}}}{{{{\left( {{\rm E}\left[ {{{\left( X \right)}^2}} \right]} \right)}^3}}}} \right]
	\end{equation*}
	The above results is then used to write the following relation: ${\rm E}\left[ {{{\left( {W({s_i}) - W({s_{i - 1}})} \right)}^6}} \right] = \mathbb{C}{\left( {{\rm E}\left[ {{{\left( {W({s_i}) - W({s_{i - 1}})} \right)}^2}} \right]} \right)^3}$. This is utilized to derive the proof as,
	\begin{equation*}
	\begin{array}{ll}
	Var[C{V_n}] &= \sum\limits_{i = 1}^n {\left( {\mathbb{C} {\rm E}{{\left[ {{{\left( {W({s_i}) - W({s_{i - 1}})} \right)}^2}} \right]}^3} } \right)} \\&
	= \mathbb{C} \sum\limits_{i = 1}^n {{{\left( {({s_i}) - ({s_{i - 1}})} \right)}^3}} \\&
	\le \mathbb{C} \underbrace {\max }_{1 \le i \le n}\left( {({s_i}) - ({s_{i - 1}})} \right)\sum\limits_{i = 1}^n {\left( {({s_i}) - ({s_{i - 1}})} \right)^2} \\&
	= \mathbb{C} \underbrace {\max }_{1 \le i \le n}\left( {({s_i}) - ({s_{i - 1}})} \right){Q_n}(id,t)\\&
	\to 0 
	\end{array}
	\end{equation*}
	Here, $\mathbb{C}$ is some constant, it measures the relation between the Hyper-Kurtosis and variance of a random process. This completes the proof. 
	In symbolic notation, without loss of generality the above cubic identities can also be expressed as, $dB_t^2dt = 0$ and $dB_t^3 = 0$. $\blacksquare$
\end{proof}

\subsection{Error process}
In order to define the additive correction, lets define an error process: ${{\mathbf{e}}_t} = {\mathbf{g}}\left( {t,{{\mathbf{X}}_t}} \right) - {\mathbf{h}}\left( {t,{{\mathbf{X}}_t}} \right)$ which follows from Eq. \ref{err_gen} such that ${\mathbf{\gamma }}\left( {t,{{\mathbf{X}}_t}} \right) = {\mathbf{f}}{\left( {t,{{\mathbf{X}}_t}} \right)^{ - 1}}{{\mathbf{e}}_t}$. To find an appropriate ${{\mathbf{e}}_t}$ for each of the Milstein schemes it is assumed that there exists an equivalent Milstein-SDE of Eqs. (\ref{explicit}), (\ref{semi_implicite_sde}) and (\ref{implicit}) that maps generic SDE (\ref{sdeg}) into explicit, semi-implicit and implicit Milstein discretization. For brevity, lets define the drift and diffusion terms in Milstein schemes in terms of ${t^*}$ such that ${t^*} = {t_{i - 1}}$ for explicit terms and ${t^*} = {t_i}$ for implicit terms. With these considerations, a Milstein-SDE is given as, 
\begin{equation}\label{sde_milstein}
\begin{array}{ll}
 d\left( {{{\mathbf{X}}_{{t_i}}}} \right) = {\mathbf{g}}\left( {{t^*},{{\mathbf{X}}_{{t^*}}}} \right)dt + {\mathbf{f}}\left( {{t^*},{{\mathbf{X}}_{{t^*}}}} \right)d{{\mathbf{W}}_t} + \frac{1}{2}{\mathbf{f}}\left( {{t^*},{{\mathbf{X}}_{{t^*}}}} \right){\mathbf{f}}{\left( {{t^*},{{\mathbf{X}}_{{t^*}}}} \right)^T}\left( {d{\mathbf{W}}_t^2 - dt} \right)
\end{array}
\end{equation}
Then, by comparing the SDEs in Eqs. (\ref{sdeg}) and (\ref{sde_milstein}) a vector valued error process is defined as, 
\begin{multline}
{{\mathbf{f}}_t}{{\mathbf{e}}_t} = \left( {{\mathbf{g}}\left( {t,{{\mathbf{X}}_t}} \right) - {\mathbf{g}}\left( {{t^*},{{\mathbf{X}}_{{t^*}}}} \right)} \right)dt + \left( {{\mathbf{f}}\left( {t,{{\mathbf{X}}_t}} \right) - {\mathbf{f}}\left( {{t^*},{{\mathbf{X}}_{{t^*}}}} \right)} \right)d{{\mathbf{W}}_t}  - \frac{1}{2}{\mathbf{f}}\left( {{t^*},{{\mathbf{X}}_{{t^*}}}} \right){\mathbf{f}}{\left( {{t^*},{{\mathbf{X}}_{{t^*}}}} \right)^T}\left( {d{\mathbf{W}}_t^2 - dt} \right)
\end{multline}
This result enables the rephrase of the Eq. (\ref{sdeg}) under the ${{({\mathcal{F}_t})}_t}$ filtered probability space $\left( {\Omega,{{({\mathcal{F}_t})}_t}} \right)$ as,
\begin{equation}\label{modsde}
\begin{array}{ll}
d{{\mathbf{X}}_t} = {\mathbf{g}}\left( {{t^*},{{\mathbf{X}}_{{t^*}}}} \right)dt + {\mathbf{f}}\left( {{t^*},{{\mathbf{X}}_{{t^*}}}} \right)\left( {{{\mathbf{e}}_t} + d{{\mathbf{W}}_t} + \frac{1}{2}{\mathbf{f}}{{\left( {{t^*},{{\mathbf{X}}_{{t^*}}}} \right)}}\left( {d{\mathbf{W}}_t^2 - dt} \right)} \right)
\end{array}
\end{equation}
After construction of proper SDE form of each of the Miltein schemes and appropriate error process, the above modified SDE can be identified as:
\begin{enumerate}
	\item Explicit Milstein:
	\begin{equation}\label{Q_brownian1}
	   d{{\bf{X}}_t} = {\bf g}\left( {{t_{i - 1}},{{\bf{X}}_{{t_{i - 1}}}}} \right)dt + {{\bf f}\left( {{t_{i - 1}},{{\bf{X}}_{{t_{i - 1}}}}} \right)}  \Bigl({{\bf e}_t} + d{{\bf W}_t}+  \cfrac{1}{2}{\bf f}\left( {{t_{i - 1}},{{\bf{X}}_{{t_{i - 1}}}}} \right) \left( {d{\bf{W}}_t^2 - dt} \right)\Bigr) 
	\end{equation}
	
	\item Semi-implicit Milstein:
	\begin{equation}\label{Q_brownian2}
	d{{\bf{X}}_t} = {\bf g}\left( {{t_i},{{\bf{X}}_{{t_i}}}} \right)dt + {{\bf f}\left( {{t_{i - 1}},{{\bf{X}}_{{t_{i - 1}}}}} \right)} \Bigl({{\bf e}_t} + d{{\bf W}_t}+  \cfrac{1}{2}{\bf f}\left( {{t_{i - 1}},{{\bf{X}}_{{t_{i - 1}}}}} \right) \left( {d{\bf{W}}_t^2 - dt} \right)\Bigr) 
	\end{equation}
	
	\item Implicit Milstein:
	\begin{equation}\label{Q_brownian3}
	d{{\mathbf{X}}_t} = {\mathbf{g}}\left( {{t_i},{{\mathbf{X}}_{{t_i}}}} \right)dt + {\mathbf{f}}\left( {{t_i},{{\mathbf{X}}_{{t_i}}}} \right) \Bigl({{\mathbf{e}}_t} + d{{\mathbf{W}}_t} +  \frac{1}{2}{\mathbf{f}}\left( {{t_i},{{\mathbf{X}}_{{t_i}}}} \right)\left( {d{\mathbf{W}}_t^2 - dt} \right) \Bigr)
	\end{equation}
\end{enumerate}
One of the challenge arises here is the separability of the terms ${\bf f}\left( {{t_{i - 1}},{{\bf{X}}_{{t_{i - 1}}}}} \right)$, which is in fact the coefficient associated with $O(h^{1.0})$ MSI in Milstein mapping of a SDE. However, in the interest of large class of dynamical vibrating oscillators appearing in Mechanical and Civil structures this is separable, which is proved in \textbf{Proposition \ref{thm1}}. 

\subsection{Weak correction on Milstein scheme}
The primary aim here is to construct a suitable Girsanov transformation for P $\to$ Q such that the additional drift ${\mathbf{f}}\left( {{t^*},{{\mathbf{X}}_{{t^*}}}} \right)\left( {{{\mathbf{e}}_t} + \frac{1}{2}{\mathbf{f}}\left( {{t^*},{{\mathbf{X}}_{{t^*}}}} \right)\left( {d{\mathbf{W}}_t^2 - dt} \right)} \right)$ gets eliminated from the modified SDE in Eq. (\ref{modsde}) and under the Q-probability space ($\Omega$,$\mathcal{F}_t$), ${{\bf{X}}_t}$ satisfies the following SDE: 
\begin{equation}\label{Q_brownian}
\begin{array}{ll}
d{{\mathbf{X}}_t} = {\mathbf{g}}\left( {{t^*},{{\mathbf{X}}_{{t^*}}}} \right)dt + {\mathbf{f}}\left( {{t^*},{{\mathbf{X}}_{{t^*}}}} \right)d{{{\mathbf{\tilde W}}}_t}
\end{array}
\end{equation}
The term ${\tilde {\bf W}}_t$ represents the zero mean Q-Wiener process independent of P, causing the transformation of ${\bf W}_t$, and defined on $\left( {\Omega ,{\mathcal{F}_t}} \right)$. The evolution of ${\bf W}_t$ is defined by the additional n-dimensional SDE:
\begin{equation}
d{\tilde {\bf W}}_t = {{\bf e}_t} + d{{\bf W}_t}+\cfrac{1}{2}{\bf f}\left( {{t^*},{{\mathbf{X}}_{{t^*}}}} \right) \left( {d{\bf{W}}_t^2 - dt} \right)
\end{equation}
where ${{\bf e}_t}$ is an additional drift of dimension-n, driving the drift field of the SDE causing the transformation between two probability measures and $\frac{1}{2}{\bf f}\left( {{t^*},{{\mathbf{X}}_{{t^*}}}} \right) \left( {d{\bf{W}}_t^2 - dt} \right)$ corresponds to the Milstein correction to the transformation (choice of ${t^*}$ will define the type of Milstein correction). To generate N-finite ensemble of Q-realizations the above equation is simulated within a Monte-Carlo (MC) framework. Thus if ${{\bf{X}}_{{t_i}}^j}$ represents the $j^{th}$ realization of ${{\bf{X}}_{{t_i}}}$ for $j=1, \ldots N$, the incremental explicit Milstein mapping is thus given as,
\begin{multline}\label{explicit_ensem}
{\mathbf{X}}_{{t_i}}^j = {\mathbf{X}}_{{t_{i - 1}}}^j + {\mathbf{g}}\left( {{t^*},{\mathbf{X}}_{_{{t^*}}}^j} \right)\Delta t +  {\mathbf{f}}\left( {{t^*},{\mathbf{X}}_{_{{t^*}}}^j} \right)\Delta \tilde {\mathbf{W}}_{{t_i}}^j + \frac{1}{2}{\mathbf{f}}\left( {{t^*},{\mathbf{X}}_{_{{t^*}}}^j} \right){\mathbf{f}}{\left( {{t^*},{\mathbf{X}}_{_{{t^*}}}^j} \right)^T}\left( {\Delta \widetilde {\mathbf{W}}_{{t_i}}^2 - \Delta t} \right)
\end{multline}
This yields a finite ensemble $\left[ {{\bf{X}}_{{t_i}}^j} \right]_{j = 1}^N$ which is the empirical distribution of ${{\bf{X}}_{{t}}}$, however, the likelihood, or weight of the particles in the distribution is decided by the Radon-Nikodym derivative, ${\bf Z}_t$. The exact solution at $t_i$ can be obtained using the empirically evaluated distribution of finite Q-ensemble of ${\bf Z}_t {{\bf{X}}_{t}}$, however as correctly indicated in \cite{beskos2005exact} the stochastic exponential ${\bf Z}_t$ causes divergence with increase in sample size due to lower weight factor. The idea here is to incorporate the ${\bf Z}_t$ as an additive term in weak sense. The weakly corrected solution of the explicit Milstein approximation of ${\bf X}_t$ can be obtained under the P probability measure by constructing a conditional expectation: ${\pi _t}({\bf{X}}) = {{\rm E}_P}\left[ {\left. {{\bf{X}}_t} \right|\mathcal{F}_t^e} \right]$, where $\mathcal{F}_t^e$ is the filtration generated by the error process ${\bf e}_t$, constructed using the sigma algebra associated with ${\bf e}_t$. Under the generated filtration the error ${\bf e}_t$ is assumed to behave as a zero mean Wiener process. Further, to ensure that the error process retains the property of a Ito-integral, the corrected solution $\pi_t({\bf{X}})$ under the measure-Q is defined as,
\begin{equation}\label{kallin}
{{\text{E}}_P}\left[ {{{\mathbf{X}}_t}|\mathcal{F}_t^e} \right] = \frac{{{{\text{E}}_Q}\left[ {\left. {{{\mathbf{Z}}_t}{{\mathbf{X}}_t}} \right|\mathcal{F}_t^e} \right]{\text{ }}}}{{{{\text{E}}_Q}\left[ {\left. {{{\mathbf{Z}}_t}} \right|\mathcal{F}_t^e} \right]}}
\end{equation}
The above expression is also called as Kallianpur-Striebel formula \cite{kallianpur1969stochastic}, where the numerator ${\sigma _t}{\text{(}}{\mathbf{X}}{\text{) = }}{{\text{E}}_Q}\left[ {\left. {{{\mathbf{Z}}_t}{{\mathbf{X}}_t}} \right|\mathcal{F}_t^e} \right]$ is often called as un-normalized estimate of ${\pi _t}({\bf{X}})$ and the denominator, ${{{\rm E}_Q}\left[ {\left. {{{\bf Z} _t}} \right|\mathcal{F}_t^e} \right]}$ acts as a normalizing constant. Alternatively, this satisfies the following equation ($d{{\bf{W}}_t} \to 0$):
\begin{equation}\label{error_evolv}
0 = -{{\bf e}_t} + {\bf \rho} d{{\tilde{\bf{W}}}_t}-\cfrac{1}{2}{\bf f}\left( {{t^*},{{\mathbf{X}}_{{t^*}}}} \right) dt 
\end{equation} 
The Eqs. (\ref{explicit_ensem}) and (\ref{error_evolv}) can be referred to as prediction and observation equations, respectively. Here, $\rho$ is an equivalent $m \times n$ diffusion matrix associated with $\mathcal{F}_t$-measurable $\tilde{{\bf W}}$ process. Comparing, Eqs. (\ref{qbrownian} ) and (\ref{error_evolv}) , one can identify, $\gamma \left( {t,{{\bf{X}}_t}} \right)= {\bf \rho}^{-1} \left({{\bf e}_t}/dt + \frac{1}{2}{\bf f}\left( {{t^*},{{\mathbf{X}}_{{t^*}}}} \right) \right)$. With ${\bf e}_t$ as zero mean martingale, the associated Radon-Nikodym derivative from Eq. (\ref{radon}) has the form:
\begin{equation}\label{radon_explicit}
    \begin{split}
        {{\mathbf{Z}}_t} = & \exp \Biggl( \int_{{t_{i - 1}}}^{{t_i}} {{{\left( {{{\mathbf{\rho }}^{ - 1}}\left( {\frac{{{{\mathbf{e}}_t}}}{{dt}} + \frac{1}{2}{{\mathbf{f}} }\left( {{t^*},{{\mathbf{X}}_{{t^*}}}} \right)} \right)} \right)}^T}} d{{{\mathbf{\tilde W}}}_s} - \\ & \frac{1}{2}{{\int_{{t_{i - 1}}}^{{t_i}} {\left\| {\int_{{t_{i - 1}}}^{{t_i}} {{{\left( {{{\mathbf{\rho }}^{ - 1}}\left( {\frac{{{{\mathbf{e}}_t}}}{{dt}} + \frac{1}{2}{{\mathbf{f}} }\left( {{t^*},{{\mathbf{X}}_{{t^*}}}} \right)} \right)} \right)}^T}} d{{{\mathbf{\tilde W}}}_s}} \right\|} }^2}ds \Biggr)
    \end{split}
\end{equation}
For generalizing the concept, one may define a bounded $\mathcal{F}_t$-measurable smooth function $\Phi ({\bf{X}})$ for ${\bf X}_t$. To obtain an evolutionary equation for $\pi_t({\bf{X}})$, one needs to apply Ito-lemma on ${{\mathbf{Z}}_t}{{\mathbf{X}}_t} \approx {{\mathbf{Z}}_t}\Phi ({{\mathbf{X}}_t})$,
\begin{equation}
\begin{array}{ll}
{{\mathbf{Z}}_t}\Phi ({{\mathbf{X}}_t}) =& \Phi ({{\mathbf{X}}_0}) + \int_0^t {{{\mathbf{Z}}_s}\Im _s^0\left( {\Phi ({{\mathbf{X}}_s})} \right)ds}  + \\ & 
\int_0^t {{{\mathbf{Z}}_s}{{\left( {\nabla \Phi \left( {{{\mathbf{X}}_s}} \right)} \right)}^T}{\mathbf{f}}\left( {s,{{\mathbf{X}}_s}} \right)d{{\mathbf{W}}_s}}  +  \int_0^t {{{\mathbf{Z}}_s}\Phi ({{\mathbf{X}}_s}){\mathbf{\gamma }}\left( {s,{{\mathbf{X}}_s}} \right)d{{{\mathbf{\tilde W}}}_s}} 
\end{array}
\end{equation}
Let ${{{\mathbf{W}}_t}}$ be a Wiener process and $\mathcal{F}_t^{\mathbf{W}} \subset {\mathcal{F}_t}$ be a $\sigma$-algebra defined as: $\mathcal{F}_t^{\mathbf{W}} = \sigma \left\{ {{{\mathbf{W}}_s}:s \leqslant t} \right\}$. If ${f\left( {{{\mathbf{X}}_t}} \right)}$ is a ${\mathcal{F}_t}$-measurable square integrable function then by Fubini's theorem: ${{\text{E}}_Q}\left[ {\int_0^t {f\left( {{{\mathbf{X}}_s}} \right)d{{\mathbf{W}}_s}|} \mathcal{F}_t^{\mathbf{W}}} \right] = \int_0^t {{{\text{E}}_Q}\left[ {f\left( {{{\mathbf{X}}_s}} \right)|\mathcal{F}_s^{\mathbf{W}}} \right]d{{\mathbf{W}}_s}} $. Then, upon taking the conditional expectation with respect to the filtration ${\mathcal{F}_t^e}$ yields,
\begin{equation}
    \begin{split}
        {{\text{E}}_Q}\left[ {\left. {{{\mathbf{Z}}_t}\Phi ({{\mathbf{X}}_t})} \right|\mathcal{F}_t^e} \right] = & {{\text{E}}_Q}\left[ {\Phi ({{\mathbf{X}}_0})} \right] + \int_0^t {{{\text{E}}_Q}\left[ {\left. {{{\mathbf{Z}}_s}\Im _s^0\left( {\Phi ({{\mathbf{X}}_s})} \right)} \right|\mathcal{F}_t^e} \right]ds}  + \\ 
        & \int_0^t {{{\text{E}}_Q}\left[ {\left. {{{\mathbf{Z}}_s}\Phi ({{\mathbf{X}}_s}){\mathbf{\gamma }}\left( {s,{{\mathbf{X}}_s}} \right)} \right|\mathcal{F}_t^e} \right]d{{{\mathbf{\tilde W}}}_s}} 
    \end{split}
\end{equation}
Recalling, the definition of un-normalized conditional estimate ${\sigma _t}{\text{(}}{\mathbf{X}}{\text{)}}$, the numerator of Eq. (\ref{kallin}) is approximated as:
\begin{equation}
{\sigma _t}\left( {\Phi ({\mathbf{X}})} \right) = {\sigma _0}\left( {\Phi ({\mathbf{X}})} \right) + \int_0^t {{\sigma _s}\left( {\Im _s^0\left( {\Phi ({\mathbf{X}})} \right)} \right)ds}  + \int_0^t {{\sigma _s}\left( {{\mathbf{\gamma }}\left( {s,{{\mathbf{X}}_s}} \right)} \right)d{{{\mathbf{\tilde W}}}_s}} 
\end{equation}
Which is the integral representation of the Zakai differential equation \cite{kallianpur2013stochastic,zakai1969optimal}. Here, $d{{{\mathbf{\tilde W}}}_t} = {\mathbf{\gamma }}\left( {t,{{\mathbf{X}}_t}} \right)dt + d{{\mathbf{W}}_t}$. The time recursive equation for normalized estimate: ${\pi _t}({\bf{X}}) = {{\rm E}_P}\left[ {\left. {{\bf{X}}_t} \right|\mathcal{F}_t^e} \right]$, is then obtained using the normalizing factor ${{{\rm E}_Q}\left[ {\left. {{{\bf Z} _t}} \right|\mathcal{F}_t^e} \right]}$ as:
\begin{equation}
    \begin{split}
        {\pi _t}\left( {\Phi ({\mathbf{X}})} \right) = & {\pi _0}\left( {\Phi ({\mathbf{X}})} \right) + \int_0^t {{\pi _s}\left( {\Im _s^0\left( {\Phi ({\mathbf{X}})} \right)} \right)ds}  + \\ & \int_0^t {\left( {{\pi _s}\left( {{\mathbf{\gamma }}\left( {s,{{\mathbf{X}}_s}} \right)\Phi ({\mathbf{X}})} \right) - {\pi _s}\left( {\Phi ({\mathbf{X}})} \right){\pi _s}\left( {{\mathbf{\gamma }}\left( {s,{{\mathbf{X}}_s}} \right)} \right)} \right)\left( {d{{{\mathbf{\tilde W}}}_s} - {\pi _s}\left( {{\mathbf{\gamma }}\left( {s,{{\mathbf{X}}_s}} \right)} \right)ds} \right)} 
    \end{split}
\end{equation}
The above equation is known as Kushner-Stratonovich or KS-equation. The weakly corrected solution of ${\bf X}_t$ is obtained as ${\pi _t}(\Phi ({\bf{X}}))$ which is in line with other KS-based stochastic filtering strategies. Then the governing evolution process of the scalar valued function ${\pi _t}(\Phi ({\bf{X}}))$ is determined following the identical filtering problems arising in nonlinear filtering theory \cite{zakai1969optimal} which is succinctly presented in Appendix \ref{filtering}. It becomes imperative at this stage to introduce the \textit{cubic variation and covariation} identities in order to derive an expression for ${\pi _t}(\Phi ({\bf{X}}))$. Following the appropriate substitutions for ${t^*}$, the process ${{\bf{\gamma }}\left( {t,{{\bf{X}}_t}} \right)}$ can be found for corresponding Milstein scheme, and utilizing the relation $d{{\mathbf{Z}}_t} = {{\mathbf{Z}}_t}{\mathbf{\gamma }}\left( {t,{{\mathbf{X}}_t}} \right)d{{{\mathbf{\tilde W}}}_t}$, the evolution expression for the conditional estimate ${\pi _t}(\Phi ({\bf{X}}))$ between the time interval $t \in \left[ {{t_{i - 1}},{t_i}} \right]$ in integral representation is given as,
\begin{equation}\label{filtration_exp}
\begin{array}{ll}
{\pi _t}\left( {\Phi ({\mathbf{X}})} \right) =& {\pi _{{t_{i - 1}}}}\left( {\Phi ({\mathbf{X}})} \right) + \int_{{t_{i - 1}}}^{{t_i}} {{\pi _s}\left( {\Im _s^0\left( {\Phi ({\mathbf{X}})} \right)} \right)ds}  + \int_{{t_{i - 1}}}^{{t_i}} {{\pi _s}\left( {\Im _s^1\left( {\Phi ({\mathbf{X}})} \right)} \right)ds}  + \\& \int_{{t_{i - 1}}}^{{t_i}} {\left( {{\pi _s}\left( {\Im _s^2\left( {\Phi \left( {\mathbf{X}} \right)} \right)} \right) - {\pi _s}\left( {{\Phi _s}\left( {\mathbf{X}} \right)} \right){\pi _s}\left( {{{\left( {{\mathbf{\gamma }}\left( {s,{{\mathbf{X}}_s}} \right)} \right)}^T}} \right)} \right)d{\psi _s}} 
\end{array}
\end{equation}
where, $d{\psi _t} = d{{{\mathbf{\tilde W}}}_t} - {\pi _t}\left( {{\mathbf{\gamma }}\left( {t,{{\mathbf{X}}_t}} \right)} \right)dt$, is called as a innovation process. The first three terms together in the approximation obtains a finite ensemble of $\Phi ({\bf{X}})$ within a MC framework using corresponding Milstein mapping and the last term is provides a weak correction to the implicit particles within a MC framework. The generators, $\Im _t^0\left( \Phi ({\bf{X}}) \right)$ is the equivalent SDE generator of Eq. (\ref{sdeg}), which is perhaps the backward Kolmogorov moment, $\Im _t^1\left( \Phi ({\bf{X}}) \right)$ is the additive order $O(\Delta t^1)$ terms from Milstein mapping and $\Im _t^2\left( \Phi ({\bf{X}}) \right)$ is an operator that accounts for the correction. Together $\Im _t^0\left( . \right)$ and $\Im _t^1\left( . \right)$ constitutes a Milstein SDE generator.

The underlying prediction equation (Eq. (\ref{explicit_ensem})), corresponding ${{\mathbf{\gamma }}\left( {t,{{\mathbf{X}}_t}} \right)}$ and the generators $\Im _t^{(.)}\left( . \right)$ in Eq. (\ref{filtration_exp}) for the Milstein schemes in section \ref{milstein_schemes} are given below,
\begin{enumerate}
	\item Girsanov corrected explicit Milstein (GCEML):\\
	For brevity lets define, ${{\mathbf{X}}_{i - 1}} = {{\mathbf{X}}_{{t_{i - 1}}}}$. Without loss of generality, through the proper substitution of ${t^*} = {t_{i-1}}$ in the drift ${\bf g}\left( {{t^*},{{\bf X}_{{t^*}}}} \right)$ and diffusion ${\bf f}\left( {{t^*},{{\bf X}_{{t^*}}}} \right)$ terms, the $j^{th}$ realization of the incremental explicit Milstein mapping follows,
	\begin{multline}
	{\mathbf{X}}_{{t_i}}^j = {\mathbf{X}}_{{t_{i - 1}}}^j + {{\mathbf{g}}_{{t_{i - 1}},{\mathbf{X}}_{_{i - 1}}^j}}\Delta t + {{\mathbf{f}}_{{t_{i - 1}},{\mathbf{X}}_{_{i - 1}}^j}}\Delta {\mathbf{\tilde W}}_{{t_i}}^j + \frac{1}{2}{{\mathbf{f}}_{{t_{i - 1}},{\mathbf{X}}_{_{i - 1}}^j}}{\mathbf{f}}_{{t_{i - 1}},{\mathbf{X}}_{_{i - 1}}^j}^T\left( {\Delta {{{\mathbf{\tilde W}}}_{{t_i}}}^2 - \Delta t} \right)
	\end{multline}
	This produces, ${\mathbf{\gamma }}\left( {t,{{\mathbf{X}}_t}} \right) = {{\mathbf{\rho }}^{ - 1}}\left( {{{\mathbf{e}}_t}/dt + \frac{1}{2}{{\mathbf{f}}_{{t_{i - 1}},{{\mathbf{X}}_{i - 1}}}}} \right)$. The generators are:
	\begin{equation}
	\begin{array}{l}
	^e\Im _t^0\left( {\Phi ({\mathbf{X}})} \right) = {\left( {{\partial _x}\Phi \left( {{{\mathbf{X}}_t}} \right)} \right)^T}{{\mathbf{g}}_{{t_{i - 1}},{{\mathbf{X}}_{i - 1}}}} + \frac{1}{2}\sum\limits_{j,k = 1}^m {\sum\limits_{l = 1}^n {{{\mathbf{f}}_{{t_{i - 1}},{{\mathbf{X}}_{i - 1}}}}{\mathbf{f}}_{{t_{i - 1}},{{\mathbf{X}}_{i - 1}}}^T} \partial _{x}^2\Phi ({\mathbf{X}})}   \\
	^e\Im _t^1\left( {\Phi ({\mathbf{X}})} \right) = \frac{1}{2}{\left( {{\partial _x}\Phi \left( {{{\mathbf{X}}_t}} \right)} \right)^T}{{\mathbf{f}}_{{t_{i - 1}},{{\mathbf{X}}_{i - 1}}}}{\mathbf{f}}_{{t_{i - 1}},{{\mathbf{X}}_{i - 1}}}^T \\
	^e\Im _t^2\left( {\Phi ({\mathbf{X}})} \right) = \Phi \left( {{{\mathbf{X}}_t}} \right){\left( {{\mathbf{\gamma }}\left( {t,{{\mathbf{X}}_t}} \right)} \right)^T} 
	\end{array} 
	\end{equation}
	
	\item Girsanov filtration on semi-implicit Milstein scheme (GCSIML):\\
	The derivation of the normalized process ${\pi _t}\left( \Phi ({\bf{X}}) \right)$ is achieved by modeling the drift term to be an implicit function of the system states i.e ${\bf g}\left( {{t^*},{{\bf X}_{{t^*}}}} \right)$ evaluated at ${t^*} = {t_i}$. As mentioned earlier a finite ensemble of Q-realization for the semi-implicit scheme can be simulated within a MC framework as,
	\begin{multline}\label{semi_implicit_ensem}
	{\mathbf{X}}_{{t_i}}^j = {\mathbf{X}}_{{t_{i - 1}}}^j + {{\mathbf{g}}_{{t_i},{\mathbf{X}}_{_i}^j}}\Delta t + {{\mathbf{f}}_{{t_{i - 1}},{\mathbf{X}}_{_{i - 1}}^j}}\Delta {\mathbf{\tilde W}}_{{t_i}}^j + \frac{1}{2}{{\mathbf{f}}_{{t_{i - 1}},{\mathbf{X}}_{_{i - 1}}^j}}{\mathbf{f}}_{{t_{i - 1}},{\mathbf{X}}_{_{i - 1}}^j}^T\left( {\Delta {{{\mathbf{\tilde W}}}_{{t_i}}}^2 - \Delta t} \right)
	\end{multline}
	The yields the process, ${\mathbf{\gamma }}\left( {t,{{\mathbf{X}}_t}} \right) = {{\mathbf{\rho }}^{ - 1}}\left( {{{\mathbf{e}}_t}/dt + \frac{1}{2}{{\mathbf{f}}_{{t_{i - 1}},{{\mathbf{X}}_{i - 1}}}}} \right)$. The generators involved in the proposed GCSIML are: 
	\begin{equation}
	\begin{array}{l}
	^s\Im _t^0\left( {\Phi ({\mathbf{X}})} \right) = {\left( {{\partial _x}\Phi \left( {{{\mathbf{X}}_t}} \right)} \right)^T}{{\mathbf{g}}_{{t_i},{{\mathbf{X}}_i}}} + \frac{1}{2}\sum\limits_{j,k = 1}^m {\sum\limits_{l = 1}^n {{{\mathbf{f}}_{{t_{i - 1}},{{\mathbf{X}}_{i - 1}}}}{\mathbf{f}}_{{t_{i - 1}},{{\mathbf{X}}_{i - 1}}}^T} \partial _{x}^2\Phi ({\mathbf{X}})}   \\
	^s\Im _t^1\left( {\Phi ({\mathbf{X}})} \right) = \frac{1}{2}{\left( {{\partial _x}\Phi \left( {{{\mathbf{X}}_t}} \right)} \right)^T}{{\mathbf{f}}_{{t_{i - 1}},{{\mathbf{X}}_{i - 1}}}}{\mathbf{f}}_{{t_{i - 1}},{{\mathbf{X}}_{i - 1}}}^T  \\
	^s\Im _t^2\left( {\Phi ({\mathbf{X}})} \right) = \Phi \left( {{{\mathbf{X}}_t}} \right){\left( {{\mathbf{\gamma }}\left( {t,{{\mathbf{X}}_t}} \right)} \right)^T} 
	\end{array}
	\end{equation}
	
	\item Girsanov filtration on implicit Milstein scheme (GCIML):\\
	In this case, both the drift ${\bf g}\left( {{t^*},{{\bf X}_{{t^*}}}} \right)$  and diffusion ${\bf f}\left( {{t^*},{{\bf X}_{{t^*}}}} \right)$ are evaluated at ${t^*} = {t_i}$. The implicit Milstein approximation to an SDE over, $t \in \left[ {{t_{i - 1}},{t_i}} \right]$, is simulated in MC framework as,
	\begin{equation}
	{\mathbf{X}}_{{t_i}}^j = {\mathbf{X}}_{{t_{i - 1}}}^j + {{\mathbf{g}}_{{t_i},{\mathbf{X}}_{_i}^j}}\Delta t + {{\mathbf{f}}_{{t_i},{\mathbf{X}}_{_i}^j}}\Delta {\mathbf{\tilde W}}_{{t_i}}^j + \frac{1}{2}{{\mathbf{f}}_{{t_i},{\mathbf{X}}_{_i}^j}}{\mathbf{f}}_{{t_i},{\mathbf{X}}_{_i}^j}^T\left( {\Delta {{{\mathbf{\tilde W}}}_{{t_i}}}^2 - \Delta t} \right)
	\end{equation}
	Identifying, ${\mathbf{\gamma }}\left( {t,{{\mathbf{X}}_t}} \right) = {{\mathbf{\rho }}^{ - 1}}\left( {{{\mathbf{e}}_t}/dt + \frac{1}{2}{{\mathbf{f}}_{{t_i},{{\mathbf{X}}_i}}}} \right)$ and modifying the corresponding drift and diffusion by ${\bf g}\left( {{t_i},{{\bf X}_{{t_i}}}} \right)$ and ${\bf f}\left( {{t_i},{{\bf X}_{{t_i}}}} \right)$, the implicit generators are obtained as,
	\begin{equation}
	\begin{array}{l}
	^i\Im _t^0\left( {\Phi ({\mathbf{X}})} \right) = {\left( {{\partial _x}\Phi \left( {{{\mathbf{X}}_t}} \right)} \right)^T}{{\mathbf{g}}_{{t_i},{{\mathbf{X}}_i}}} + \frac{1}{2}\sum\limits_{j,k = 1}^m {\sum\limits_{l = 1}^n {{{\mathbf{f}}_{{t_i},{{\mathbf{X}}_i}}}{\mathbf{f}}_{{t_i},{{\mathbf{X}}_i}}^T} \partial _x^2\Phi ({\mathbf{X}})}   \\
	^i\Im _t^1\left( {\Phi ({\mathbf{X}})} \right) = \frac{1}{2}{\left( {{\partial _x}\Phi \left( {{{\mathbf{X}}_t}} \right)} \right)^T}{{\mathbf{f}}_{{t_i},{{\mathbf{X}}_i}}}{\mathbf{f}}_{{t_i},{{\mathbf{X}}_i}}^T  \\
	^i\Im _t^2\left( {\Phi ({\mathbf{X}})} \right) = \Phi \left( {{{\mathbf{X}}_t}} \right){\left( {{\mathbf{\gamma }}\left( {t,{{\mathbf{X}}_t}} \right)} \right)^T} 
	\end{array}
	\end{equation}
\end{enumerate}


\section{Numerical illustrations:}
Three nonlinear oscillators are considered in the present study for the illustration of the applicability of the proposed weakly corrected schemes. These oscillators are widely used in order to characterize the non-linearity associated with various natural mechanical and structural systems. Successful capturing of the actual dynamics under given initial and forcing function is the prime motivation here.

\subsection{Illustration 1: Duffing-Van der pol oscillator (DV)}
The dynamics of a Duffing-Van der pol oscillator is expressed by the following equation \cite{tripura2020ito}:
\begin{equation}\label{eq_sdof_composite}
\ddot X(t) + \dot X(t) -  (\alpha - {X^2}(t))X(t) = {\sigma }X(t)\dot W(t)
\end{equation}
With a statespace of $X_1=X$, and $X_2= \dot{X}$, the first order incremental Ito-diffusion equations for the Eq. (\ref{eq_sdof_composite}) can be written as:
\begin{equation}\label{diffusion_composite}
d\left[ {\begin{array}{*{20}{c}}
	{{X_1}(t)} \\ 
	{{X_2}(t)} 
	\end{array}} \right] = \left[ {\begin{array}{*{20}{c}}
	{{X_2}(t)} \\ 
	{\left( {\alpha  - X_1^2(t)} \right){X_1}(t) - {X_2}(t)} 
	\end{array}} \right]dt + \left[ {\begin{array}{*{20}{c}}
	0 \\ 
	{\sigma {X_1}(t)} 
	\end{array}} \right]dW(t)
\end{equation}
The Milstein mappings for the statespace equation within time interval, $\Delta t=(t_i - t_{i-1})$, can be formulated using Eqs. \ref{explicit}, \ref{semi_implicite_sde} and \ref{implicit} as follows:
\begin{equation}\label{SDVPcomp}
\begin{array}{l}
\text{Explicit}: \\
 {\begin{array}{*{20}{l}}
	{{X_1}\left( {{t_{i + 1}}} \right) = {X_1}\left( {{t_i}} \right) + {X_2}\left( {{t_i}} \right)\Delta t} \\ 
	{{X_2}\left( {{t_{i + 1}}} \right) = {X_2}\left( {{t_i}} \right) + \left( {\left( {\alpha  - X_1^2\left( {{t_i}} \right)} \right){X_1}\left( {{t_i}} \right) - {X_2}\left( {{t_i}} \right)} \right)\Delta t + \sigma {X_1}\left( {{t_i}} \right)\Delta {W_t}} 
	\end{array}} \\
\text{Semi-implicit}: \\
 {\begin{array}{*{20}{l}}
	{{X_1}\left( {{t_{i + 1}}} \right) = {X_1}\left( {{t_i}} \right) + {X_2}\left( {{t_{i + 1}}} \right)\Delta t} \\ 
	{{X_2}\left( {{t_{i + 1}}} \right) = {X_2}\left( {{t_i}} \right) + \left( {\left( {\alpha  - X_1^2\left( {{t_{i + 1}}} \right)} \right){X_1}\left( {{t_{i + 1}}} \right) - {X_2}\left( {{t_{i + 1}}} \right)} \right)\Delta t + \sigma {X_1}\left( {{t_i}} \right)\Delta {W_t}} 
	\end{array}}\\
\text{Implicit}: \\
 {\begin{array}{*{20}{l}}
	{{X_1}\left( {{t_{i + 1}}} \right) = {X_1}\left( {{t_i}} \right) + {X_2}\left( {{t_{i + 1}}} \right)\Delta t} \\ 
	{{X_2}\left( {{t_{i + 1}}} \right) = {X_2}\left( {{t_i}} \right) + \left( {\left( {\alpha  - X_1^2\left( {{t_{i + 1}}} \right)} \right){X_1}\left( {{t_{i + 1}}} \right) - {X_2}\left( {{t_{i + 1}}} \right)} \right)\Delta t + \sigma {X_1}\left( {{t_{i + 1}}} \right)\Delta {W_t}} 
	\end{array}} 
\end{array}
\end{equation}
The innovation term on the approximation of the proposed weakly corrected schemes can be incorporated by properly estimating the term ${{\bf{\gamma }}\left( {s,{{\bf{X}}_s}} \right)}$, which is for the undertaken oscillator takes the form,
\begin{equation}
{\mathbf{\gamma }}\left( {s,{{\mathbf{X}}_s}} \right) = {{\mathbf{\rho }}^{ - 1}}\left( {\frac{{{{\mathbf{e}}_t}}}{{dt}} + \frac{{\sigma {X_1}({t^*})}}{2}} \right)
\end{equation} 
\begin{figure*}[ht]
	\centering
	\includegraphics[width=\textwidth]{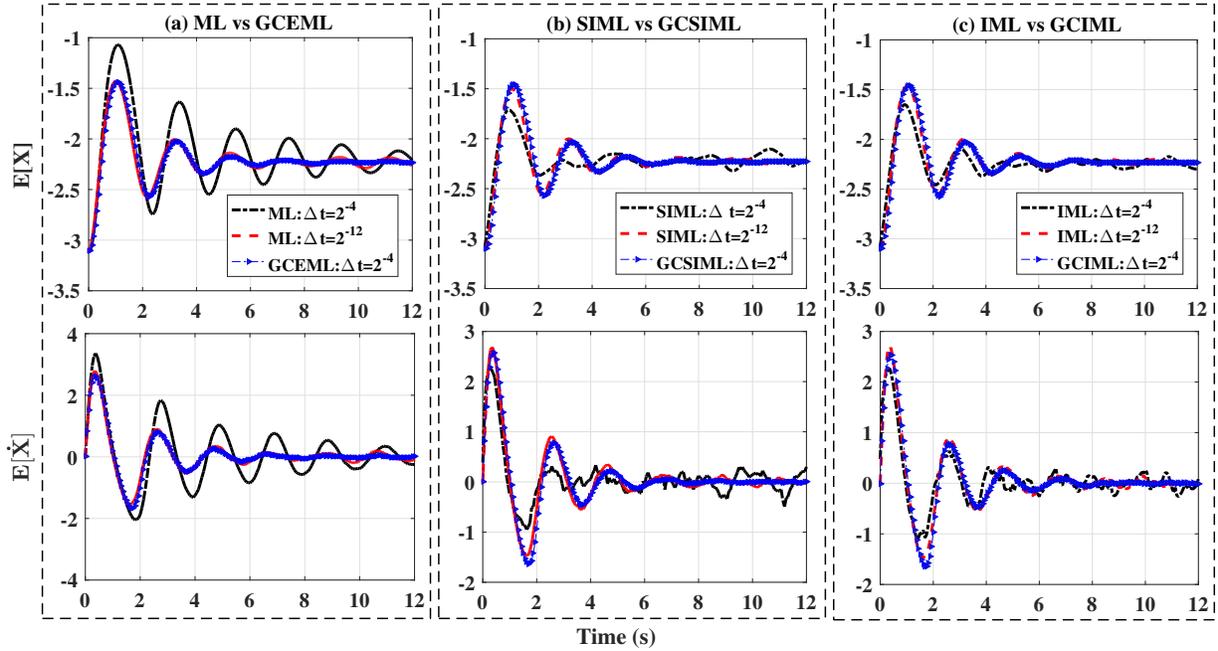}
	\caption{Duffing-Van der pol oscillator (Eq. (\ref{eq_sdof_composite})): ensemble mean displacement ({\rm E[X]}) and ensemble mean velocity ({\rm E[$\dot{X}$]}) response for $\alpha$ = 5.0 and $\sigma$ = 0.2.}
	\label{figdisvel}
\end{figure*}
Where, $\sigma {X_1}({t^*}) = \sigma {X_1}({t_i})$ for explicit and semi-implicit and $\sigma {X_1}({t^*}) = \sigma {X_1}({t_{i + 1}})$ for implicit Milstein scheme. The associated Radon-Nikodym derivative then follows Eq. (\ref{radon}). The error process ${{\bf e}_t}$ associated with the Radon-Nikodym derivative can be identified as: 
\begin{equation}
\begin{array}{ll}
\multicolumn{2}{l}{{\text{Explicit Milstein,}}} \\ 
{e_t} =& \left\{ {\left( {\left( {\alpha  - X_1^2(t)} \right){X_1}(t) - {X_2}(t)} \right) - \left( {\left( {\alpha  - X_1^2\left( {{t_i}} \right)} \right){X_1}\left( {{t_i}} \right) - {X_2}\left( {{t_i}} \right)} \right)} \right\}dt \\& +
\left\{ {\sigma X(t) - \sigma X\left( {{t_i}} \right)} \right\}d{W_t}  - \frac{1}{2}{\left( {\sigma {X_1}\left( {{t_i}} \right)} \right)^2}\left( {\Delta {{{\tilde W}}_{t}}^2 - \Delta t} \right) \\
\multicolumn{2}{l}{{\text{Semi-implicit Milstein,}}} \\
{e_t} =& \left\{ {\left( {\left( {\alpha  - X_1^2(t)} \right){X_1}(t) - {X_2}(t)} \right) - \left( {\left( {\alpha  - X_1^2\left( {{t_{i + 1}}} \right)} \right){X_1}\left( {{t_{i + 1}}} \right) - {X_2}\left( {{t_{i + 1}}} \right)} \right)} \right\}dt \\& + \left\{ {\sigma X(t) - \sigma X\left( {{t_i}} \right)} \right\}d{W_t} - \frac{1}{2}{\left( {\sigma {X_1}\left( {{t_i}} \right)} \right)^2}\left( {\Delta {{{\tilde W}}_{t}}^2 - \Delta t} \right)\\
\multicolumn{2}{l}{{\text{Implicit Milstein,}}} \\
{e_t} =& \left\{ {\left( {(\alpha  - X_1^2(t)){X_1}(t) - {X_2}(t)} \right) - \left( {\left( {\alpha  - X_1^2\left( {{t_{i + 1}}} \right)} \right){X_1}\left( {{t_{i + 1}}} \right) - {X_2}\left( {{t_{i + 1}}} \right)} \right)} \right\}dt \\& + \left\{ {\sigma X(t) - \sigma X\left( {{t_{i + 1}}} \right)} \right\}d{W_t} - \frac{1}{2}{\left( {\sigma {X_1}\left( {{t_{i + 1}}} \right)} \right)^2}\left( {\Delta {{{\tilde W}}_{t}}^2 - \Delta t} \right)
\end{array}
\end{equation}
With the initial conditions as $(X_0,\dot{X}_0)$=(-3.1,0), the system is simulated for $\alpha$=5.0 and $\sigma$=0.2 using explicit, semi-implicit, and implicit Milstein schemes for $\Delta t$=$2^{-4}$s. The weak correction through Girsanov transform is performed after simulating an ensemble of 200 realizations. The ensemble average solutions of system states are compared with a reference solution, in the absence of true solution. The reference solution is obtained though explicit Milstein scheme using a $\Delta t$=$2^{-12}$s. 
\begin{figure*}[!htb]
	\centering
	\includegraphics[width=\textwidth]{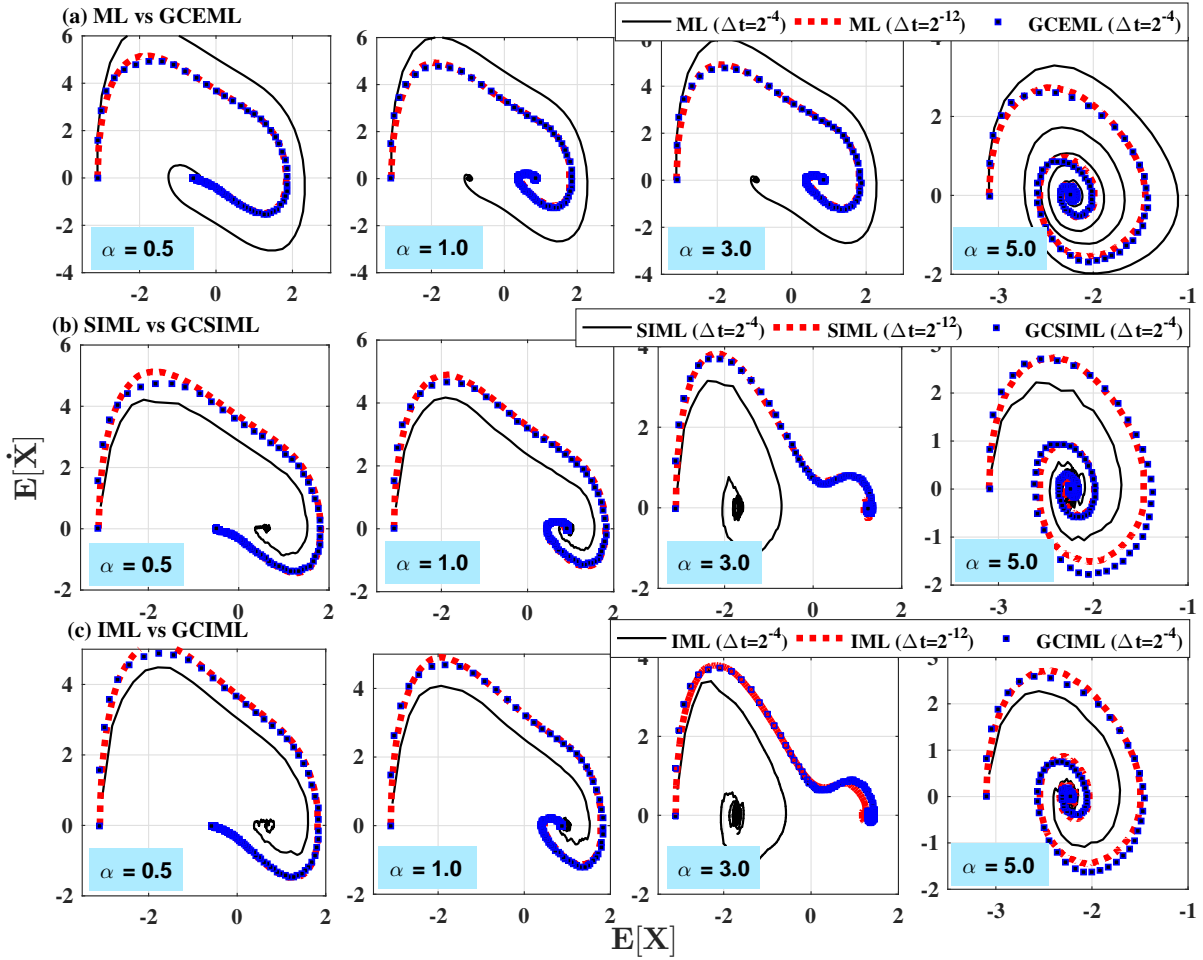}
	\caption{Duffing-Van der pol oscillator (Eq. (\ref{eq_sdof_composite})): ensemble mean stability of the response under additive noise strength $\sigma$ = 0.2 for $\alpha$= 0.5, 1.0, 3.0, and 5.0.}
	\label{figpp}
\end{figure*}
From Fig. \ref{figdisvel}, it can be observed that the Girsanov correction over the solution obtained using explicit, semi-implicit, and implicit Milstein schemes for large time step $\Delta t$=$2^{-4}$s indistinguishably evolves with time in contrast to the reference solution, whereas, the corresponding solutions without correction for $\Delta t$=$2^{-4}$s tends to suffer convergence issues. Further, in order to asses the stability of the system, the ensemble mean displacements are plotted against ensemble mean velocity trajectories for varying levels of $\alpha$ in Fig. \ref{figpp}. The steady state solutions of the system is given by $X_1(t)=\pm \sqrt{\alpha}$, and $X_2(t)=$0, however, convergence to the original steady state under small $\sigma$ entirely depends on the initial conditions. For the stated initial conditions, the phasespace trajectories obtained from the Girsanov corrected solutions always remain close to the reference solution and finally achieving the original steady state. On the contrary, although the solutions of Milstein schemes without correction lie in the vicinity of the reference solution at the beginning, however, fails to converge to the original steady state.

\subsection{Illustration 2: Duffing Holmes oscillator (DH)}
A Duffing-Holmes (DH) Oscillator excited under the presence of combined deterministic and additive noise is considered. For $\epsilon_3$ and $\epsilon_4$=0, the system has two stable steady state solutions at
($X$,$\dot{X}$) = ($\pm$1,0). This oscillator is useful for understanding the nonlinear dynamics of a periodically forced buckled beam. The governing equation of the oscillator is,
\begin{equation}\label{DH}
\begin{array}{ll}
\ddot X(t) + 2\pi {\varepsilon _1}\dot X(t) + 4{\pi ^2}{\varepsilon _2}X(t)\left( { - 1 + {X(t)^2}} \right) = 4{\pi ^2}{\varepsilon _3}\cos \left( {2\pi t} \right) + 4{\pi ^2}{\varepsilon _4}{{\dot W}_1}\left( t \right)
\end{array}
\end{equation}
For a transformation of the form, ${Y_1} = X$ and ${Y_2} = \dot X$ the statespace model of the oscillator in the form of incremental Ito-diffusion SDEs can be written as,
\begin{equation}
d\left[ {\begin{array}{*{20}{c}}
	{{X_1}} \\ 
	{{X_2}} 
	\end{array}} \right] = \left[ {\begin{array}{*{20}{c}}
	{{Y_2}} \\ 
	{ - 2\pi {\varepsilon _1}{Y_2} - 4{\pi ^2}{\varepsilon _2}{Y_1}\left( { - 1 + {Y_1}^2} \right) + 4{\pi ^2}{\varepsilon _3}\cos \left( {2\pi t} \right)} 
	\end{array}} \right]dt + \left[ {\begin{array}{*{20}{c}}
	0 \\ 
	{4{\pi ^2}{\varepsilon _4}} 
	\end{array}} \right]d{W_t}
\end{equation}

\begin{figure}[ht]
	\centering
	\includegraphics[width= \textwidth]{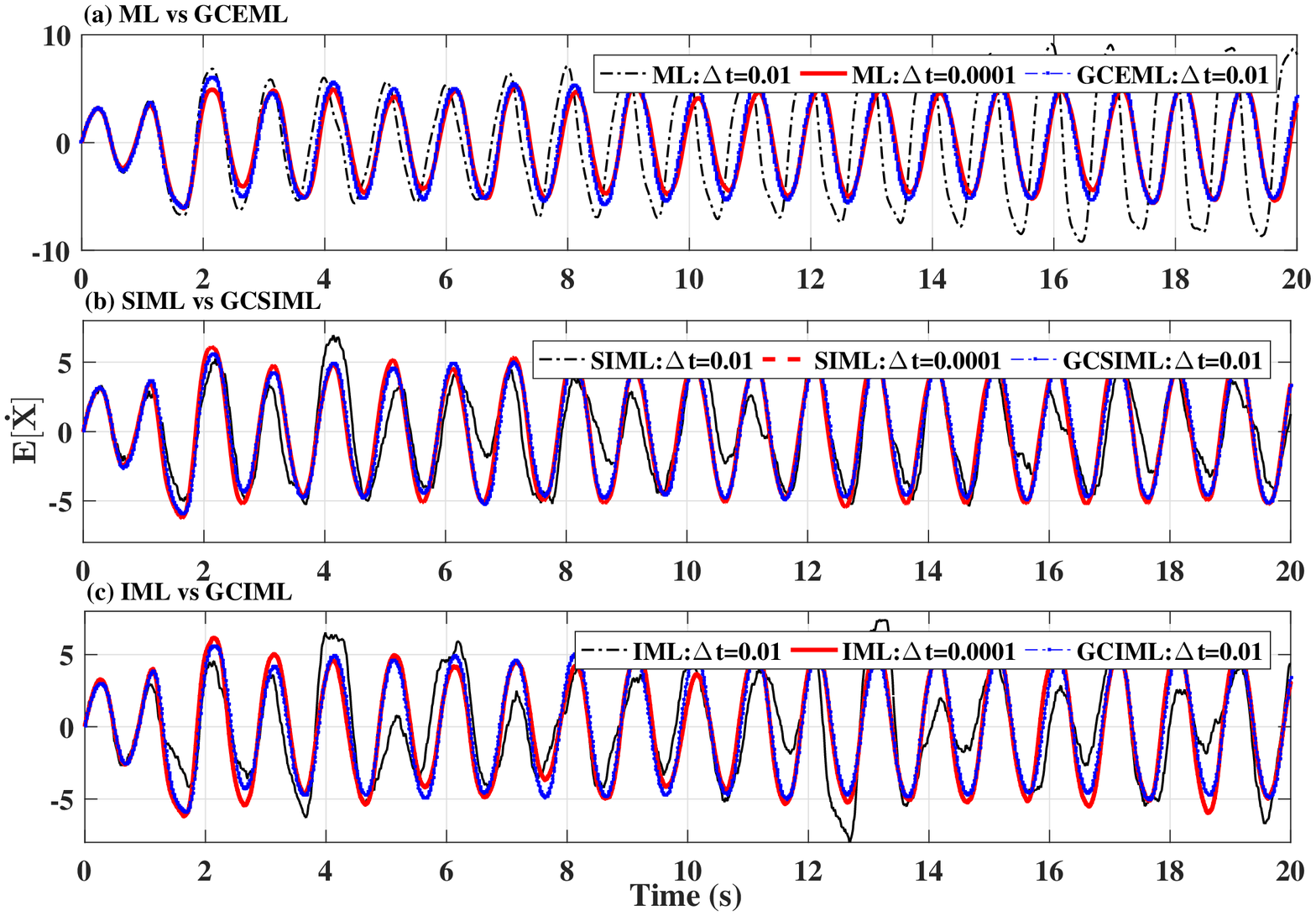}
	\caption{DH oscillator (Eq. (\ref{DH})): ensemble mean velocity ({\rm E[$\dot{X}$]}) time history $\epsilon_1$=0.25, $\epsilon_2$=0.5, $\epsilon_3$=0.5 and $\epsilon_4$=0.05.}
	\label{figveldh}
\end{figure}

\begin{figure}[!htb]
	\centering
	\includegraphics[width= \textwidth]{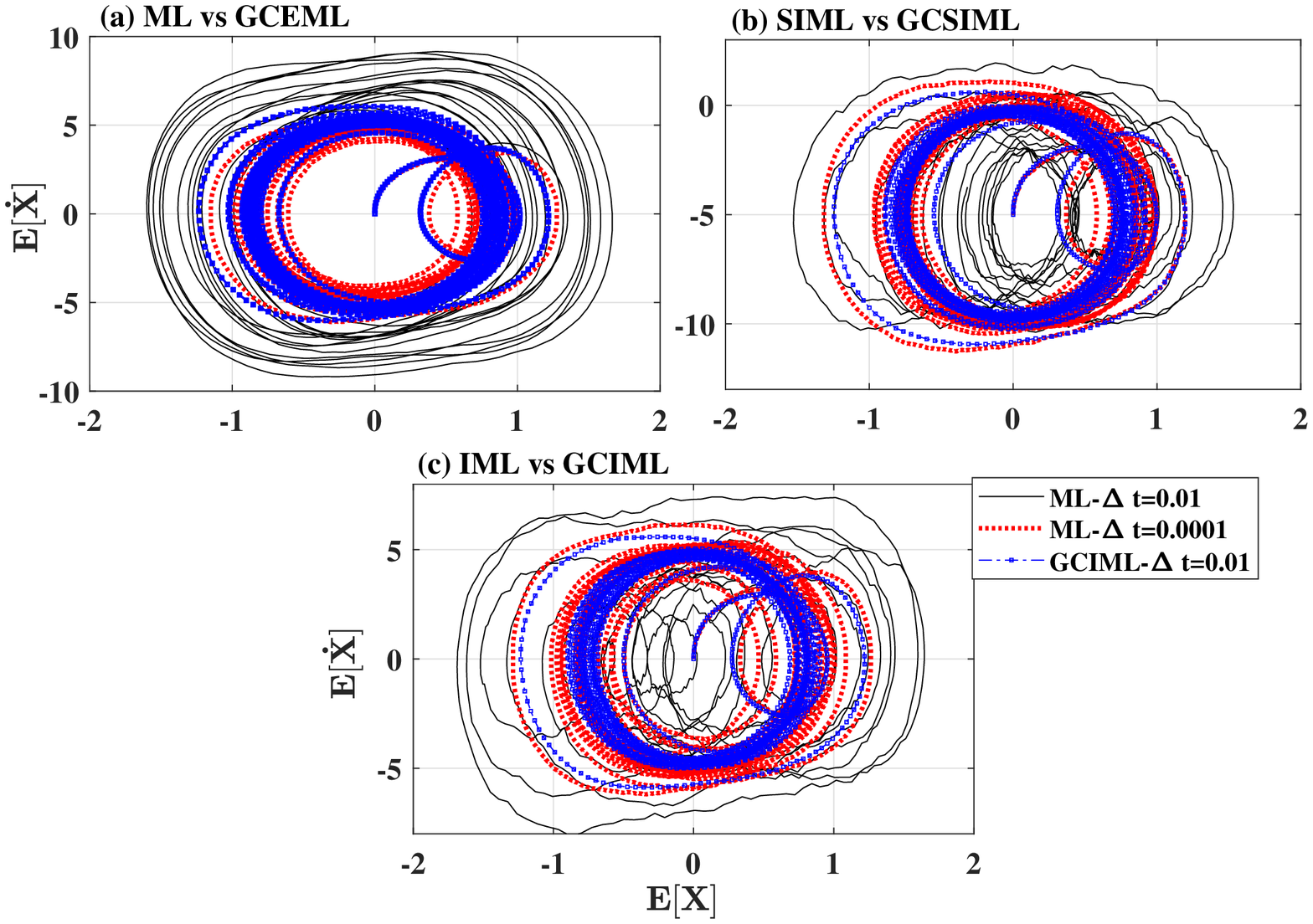}
	\caption{GCEML - DH oscillator (Eq. (\ref{DH})): ensemble mean displacement ({\rm E[X]}) vs ensemble mean velocity ({\rm E[$\dot{X}$]}).}
	\label{figppdhexp}
\end{figure}

Noting that, ${\Im ^1}{f_1}\left( {t,{{\bf{X}}_t}} \right) = {\Im ^1}{f_2}\left( {t,{{\bf{X}}_t}} \right) = 0$, the Milstein mappings for the diffusions are obtained as,
\begin{equation}
\label{DHcomp}
\begin{array}{l}
\text{Explicit}: \\
{\begin{array}{ll}
	{Y_1}\left( {{t_{i + 1}}} \right) =& {Y_1}\left( {{t_i}} \right) + {Y_2}\left( {{t_i}} \right)\Delta t\\
	{Y_2}\left( {{t_{i + 1}}} \right) =& {Y_2}\left( {{t_i}} \right) + \Bigl[ - 2\pi {\varepsilon _1}{Y_2}\left( {{t_i}} \right) - 4{\pi ^2}{\varepsilon _2}{Y_1}\left( {{t_i}} \right)\left( { - 1 + {Y_1}^2\left( {{t_i}} \right)} \right) + \\& 4{\pi ^2}{\varepsilon _3}\cos \left( {2\pi t} \right)\Bigr]\Delta t  + 4{\pi ^2}{\varepsilon _4}\Delta {W_n}
	\end{array}} \\
\text{Semi-implicit}: \\
{\begin{array}{ll}
	{Y_1}\left( {{t_{i + 1}}} \right) =& {Y_1}\left( {{t_i}} \right) + {Y_2}\left( {{t_{i + 1}}} \right)\Delta t\\
	{Y_2}\left( {{t_{i + 1}}} \right) =& {Y_2}\left( {{t_i}} \right) + \Bigl[ - 2\pi {\varepsilon _1}{Y_2}\left( {{t_{i + 1}}} \right) - 4{\pi ^2}{\varepsilon _2}{Y_1}\left( {{t_{i + 1}}} \right)\left( { - 1 + {Y_1}^2\left( {{t_{i + 1}}} \right)} \right) + \\& 4{\pi ^2}{\varepsilon _3}\cos \left( {2\pi t} \right) \Bigr]\Delta t  + 4{\pi ^2}{\varepsilon _4}\Delta {W_n}
	\end{array}} \\
\text{Implicit}: \\
{\begin{array}{ll}
	{Y_1}\left( {{t_{i + 1}}} \right) =& {Y_1}\left( {{t_i}} \right) + {Y_2}\left( {{t_{i + 1}}} \right)\Delta t\\
	{Y_2}\left( {{t_{i + 1}}} \right) =& {Y_2}\left( {{t_i}} \right) + \Bigl[ - 2\pi {\varepsilon _1}{Y_2}\left( {{t_{i + 1}}} \right) - 4{\pi ^2}{\varepsilon _2}{Y_1}\left( {{t_{i + 1}}} \right)\left( { - 1 + {Y_1}^2\left( {{t_{i + 1}}} \right)} \right) + \\& 4{\pi ^2}{\varepsilon _3}\cos \left( {2\pi t} \right) \Bigr]\Delta t  + 4{\pi ^2}{\varepsilon _4}\Delta {W_n}
	\end{array}} 
\end{array}
\end{equation}
\begin{figure}[!htb]
	\centering
	\includegraphics[width= \textwidth]{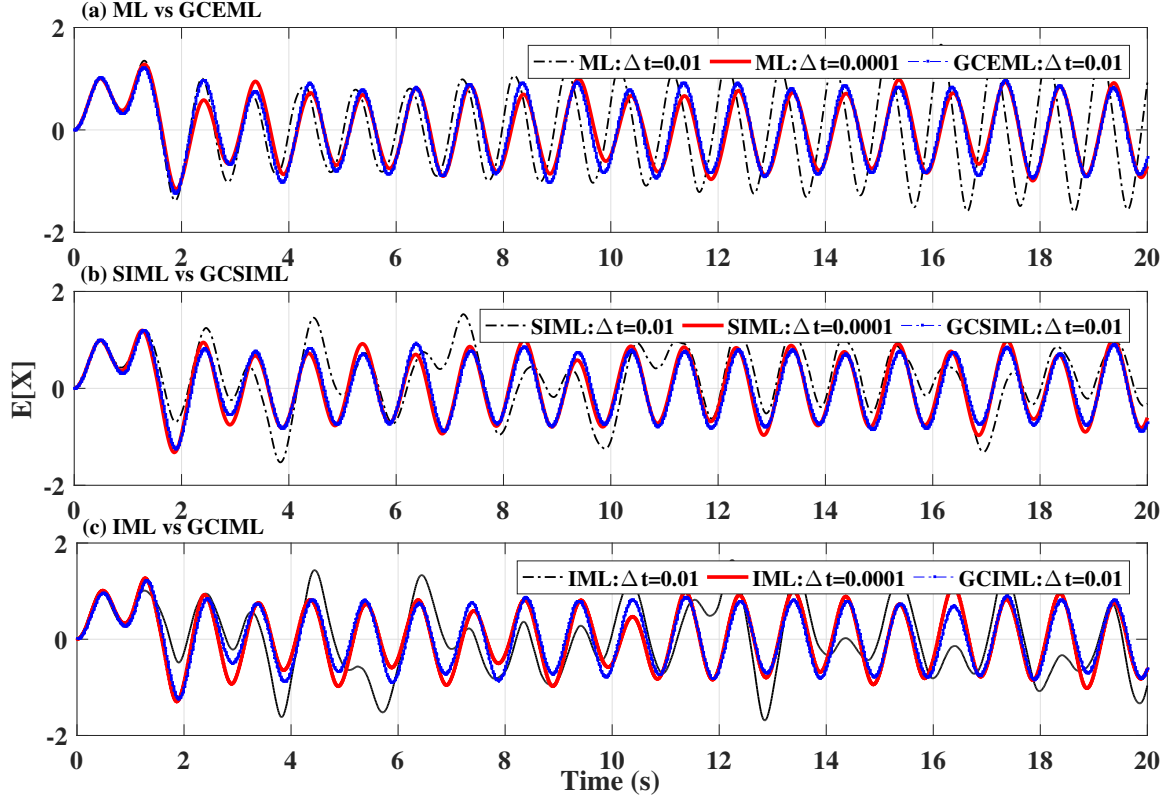}
	\caption{DH oscillator (Eq. (\ref{DH})): ensemble mean displacement ({\rm E[X]}) time history for $\epsilon_1$=0.25, $\epsilon_2$=0.5, $\epsilon_3$=0.5 and $\epsilon_4$=0.05.}
	\label{figdisdh}
\end{figure}
Following the lines of development of the proposed schemes, the error associated with $\gamma \left( {s,{{\bf{X}}_s}} \right)$ is noted as follows:
\begin{multline}
{e_t} = \bigl\{ \left( { - 2\pi {\varepsilon _1}{Y_2}\left( t \right) - 4{\pi ^2}{\varepsilon _2}{Y_1}\left( t \right)\left( { - 1 + {Y_1}^2\left( t \right)} \right)} \right) - \\ \left( { - 2\pi {\varepsilon _1}{Y_2}\left( {{t^*}} \right) - 4{\pi ^2}{\varepsilon _2}{Y_1}\left( {{t^*}} \right)\left( { - 1 + {Y_1}^2\left( {{t^*}} \right)} \right)} \right) \bigr\}dt - \frac{1}{2}{\left( {4{\pi ^2}{\varepsilon _4}} \right)^2}\left( {\Delta {{\tilde W}_t}^2 - \Delta t} \right)
\end{multline}
As discussed earlier, substitution for ${t^*} = {t_{i}}$ and ${t^*} = {t_{i + 1}}$ in above expression yields error process ${e_t}$ corresponding to explicit and semi-implicit in the former case and to implicit Milstein scheme for the later case. The system is simulated for $\epsilon_1$=0.25, $\epsilon_2$=0.5 and $\epsilon_3$=0.5 within a Monte Carlo (MC) setup using the three proposed weakly Girsanov corrected Milstein schemes. The oscillator is excited using zero mean Gaussian white noise ($\mathcal{N}~[0,1]$) with a intensity of $\epsilon_4$=0.05. The innovation term for the Girsanov correction is implemented by simulating 200 ensembles from the corresponding uncorrected Milstein schemes. The corrected solutions are compared with a reference solution that is obtained using explicit Milstein scheme with a $\Delta t$=0.0001s. The system states are portrayed in Figs. \ref{figdisdh} and \ref{figveldh}, while Fig. \ref{figppdhexp} illustrates the ensemble mean stability of the system under ($Y_1$(0),$Y_2$(0))=(0,0).

It clear that the Girsanov corrected weak solutions of GCEML matches with the reference solution with a coarser $\Delta t$=0.01s as compared to ML with same time sampling (Figs. \ref{figdisdh} and \ref{figveldh}). Further, it is seen that the SIML and IML schemes are better in approximating the reference solution than ML (explicit scheme), which suffers from numerical steadiness issues at intermediate instants. This problem is absent in corresponding proposed GCSIML and GCIML approaches. In Fig. \ref{figppdhexp}, it is observed that the solution of the ML, SIML and IML schemes does not converge to the stable orbits for the given zero initial conditions. In case of the solutions simulated using GCEML, GCSIML and GCIML schemes, the response histories always stays with the reference solution in the final steady state phasespace orbit, devoid of coarser $\Delta t$. 

\subsection{Illustration 3: Ring-Type MEMS Gyroscopes Subjected to Stochastic Angular Speed Fluctuation}
A two degrees-of-freedom ring-type microelectromechanical systems (MEMS) gyroscope is investigated assuming the effect of stochastic fluctuations in angular velocity \cite{asokanthan2017stability}. These angular rate gyroscopes are used either as a stand-alone unit or as part of an inertial measurement unit (IMU) in many applications such as in automotive traction control systems, ride stabilization, and rollover detection, in consumer electronic applications like stabilization of digital video cameras, in military applications like guidance of missiles and platform stabilization, in aerospace and marine and many more. This study is useful because during the service life it is subjected to external forces such as fabrication-deployment-and-operation impacts, vibratory excitations resulting from the operating environment.

The coordinate vector $Q$ for the ring vibration is assumed as, $Q=[q_1, q_2]^T$. During functioning of a gyroscope an angular shift between $0^\circ - 45^\circ$ is realized, under the effect the coordinate $q_1$ is considered to represent the excitation and the coordinate $q_2$ is assumed to be associated with the angular rate measurement. Under the defined coordinates, the governing equation of motion for the gyroscope is given as \cite{asokanthan2017stability}:
\begin{multline}\label{gyroscope}
\left[ {\begin{array}{*{20}{c}}
	1&0\\
	0&{1 + \delta m}
	\end{array}} \right]\ddot Q(t) +  \left[ {\begin{array}{*{20}{c}}
	{2\xi {\omega _{01}}}&{ - 2\Omega \gamma }\\
	{2\Omega \gamma }&{2\xi {\omega _{02}}}
	\end{array}} \right]\dot Q(t) + 
 \left[ {\begin{array}{*{20}{c}}
	{{\kappa _1} + {\kappa _1}{\Omega ^2}}&{ - \dot \Omega \gamma }\\
	{\dot \Omega \gamma }&{{\kappa _1} + {\kappa _1}{\Omega ^2}}
	\end{array}} \right]Q(t) = \left\{ {\begin{array}{*{20}{c}}
	{P\cos \omega t}\\
	0
	\end{array}} \right\}
\end{multline}
If constants $r$, $E$ and $\rho$ refers to ring radius, Young’s
modulus of elasticity and ring density, respectively, while $\tilde{h}$ is radial and $\tilde{b}$ is the axial thickness then, the constants are evaluated as,
\begin{equation}
\begin{array}{ll}
\gamma  = \cfrac{{\tilde b + 4\tilde a}}{{2(\tilde a + \tilde b)}}; \quad {\kappa _1} = \cfrac{{\tilde b\tilde c + 4{{\tilde a}^2}}}{{\rho A(\tilde a + \tilde b)}}; \quad {\kappa _2} = \cfrac{{4\left( {\tilde b + \tilde c - 4\tilde a}; \right)}}{{(\tilde a + \tilde b)}} - \cfrac{{4\left( {\tilde b\tilde c - 4\tilde a} \right)}}{{{{(\tilde a + \tilde b)}^2}}};\\
\tilde a = \cfrac{{4EI}}{{{r^4}}} + \cfrac{{EA}}{{{r^2}}}; \quad \tilde b = \cfrac{{4EI}}{{{r^4}}} + \cfrac{{4EA}}{{{r^2}}}; \quad \tilde c = \cfrac{{16EI}}{{{r^4}}} + \cfrac{{EA}}{{{r^2}}};
\quad A = \tilde b\tilde h; \quad I = \cfrac{{\tilde b{{\tilde h}^3}}}{{12}}
\end{array}
\end{equation}
The term $\delta m$ represents the mass mismatch of the ring, however, is assumed as $0$ to include uniformly distributed mass along the circumference of the ring. $\omega_{01}$ and $\omega_{02}$ are the two undamped system natural frequencies which depends on the input angular velocity $\Omega$, however, at typical low $\Omega= 2\pi$rad/s, takes nearly identical values. At higher angular rates the magnitude of ${\dot \Omega \gamma }$ becomes negligible as compared to $2\Omega \gamma$, and thus the constant angular rate ${\dot{\Omega}}$ is assumed as 0. To form the first-order Ito-diffusions or SDEs four state variables are introduced: $q_1=X_1$, $\dot{q}_ 1=X_2$, $q_2=X_3$ and $\dot{q}_2=X_4$. The SDEs are given as,
\begin{equation}
d\left[ {\begin{array}{*{20}{c}}
	{{X_1}}\\
	{{X_2}}\\
	{{X_3}}\\
	{{X_4}}
	\end{array}} \right] = \left[ {\begin{array}{*{20}{c}}
	{{X_2}}\\
	{ - \left( {{\kappa _1} + {\kappa _1}{\Omega ^2}} \right){X_1} - 2\xi {\omega _{01}}{X_2} + 2\Omega \gamma {X_4} + \tilde{P}}\\
	{{X_4}}\\
	{ - 2\Omega \gamma {X_2} - \left( {{\kappa _1} + {\kappa _1}{\Omega ^2}} \right){X_3} - 2\xi {\omega _{02}}{X_4}}
	\end{array}} \right]dt
\end{equation}
where, $\tilde{P}=P\cos \omega t$. In the presence of external
noise resulting from environment factors and the nature of operation can 
excite the system at an arbitrary frequency range depending on the source. In order to account for the random fluctuations in the dynamic behavior of gyroscope, the angular velocity is written in terms of nominal angular velocity ${\Omega _0}$ under additive stochastic noise, where the noise simulated using Gaussion random variables. The random fluctuation in $\Omega_0$ is expressed as, $\Omega  = {\Omega _0} + {\mu _0}\dot W(t)$. Here, ${\mu _0}$ is the strength of the noise. In the expansion ${\Omega ^2} = \Omega _0^2 + 2{\Omega _0}{\mu _0}\dot W(t) + {( {{\mu _0}\dot W(t)} )^2}$, ${( {{\mu _0}\dot W(t)})^2}$ is neglected since ${\mu _0}\dot W(t)<<1$. For more appropriate representation of the fluctuations in ${\mu_0}$, a noise intensity ration $\mu$ is introduced,
\begin{equation}
\mu  = \cfrac{{{\mu _0}(\dot W{{(t)})_{\max }}}}{{{{({\Omega _0})}_{\max }}}}
\end{equation}
The SDEs after incorporating the fluctuation in the $\Omega_0$, rephrased as,
\begin{equation}
\begin{array}{ll}
d \left[ {\begin{array}{*{20}{c}}
	{{X_1}}\\
	{{X_2}}\\
	{{X_3}}\\
	{{X_4}}
	\end{array}} \right] =& \left[ {\begin{array}{*{20}{c}}
	{{X_2}}\\
	{ - \left( {{\kappa _1} + {\kappa _1}{\Omega ^2}} \right){X_1} - 2\xi {\omega _{01}}{X_2} + 2\Omega \gamma {X_4} + \tilde{P}}\\
	{{X_4}}\\
	{ - 2\Omega \gamma {X_2} - \left( {{\kappa _1} + {\kappa _1}{\Omega ^2}} \right){X_3} - 2\xi {\omega _{02}}{X_4}}
	\end{array}} \right]dt  + \\& \left[ {\begin{array}{*{20}{c}}
	0\\
	{ - 2{\Omega _0}{\kappa _2}{\mu _0}{X_1} + 2{\mu _0}\gamma {X_4}}\\
	0\\
	{ - 2{\mu _0}\gamma {X_2} - 2{\Omega _0}{\kappa _2}{\mu _0}{X_3}}
	\end{array}} \right]dW(t)
\end{array}
\end{equation}

\begin{figure}[!htb]
	\centering
	\includegraphics[width= \textwidth]{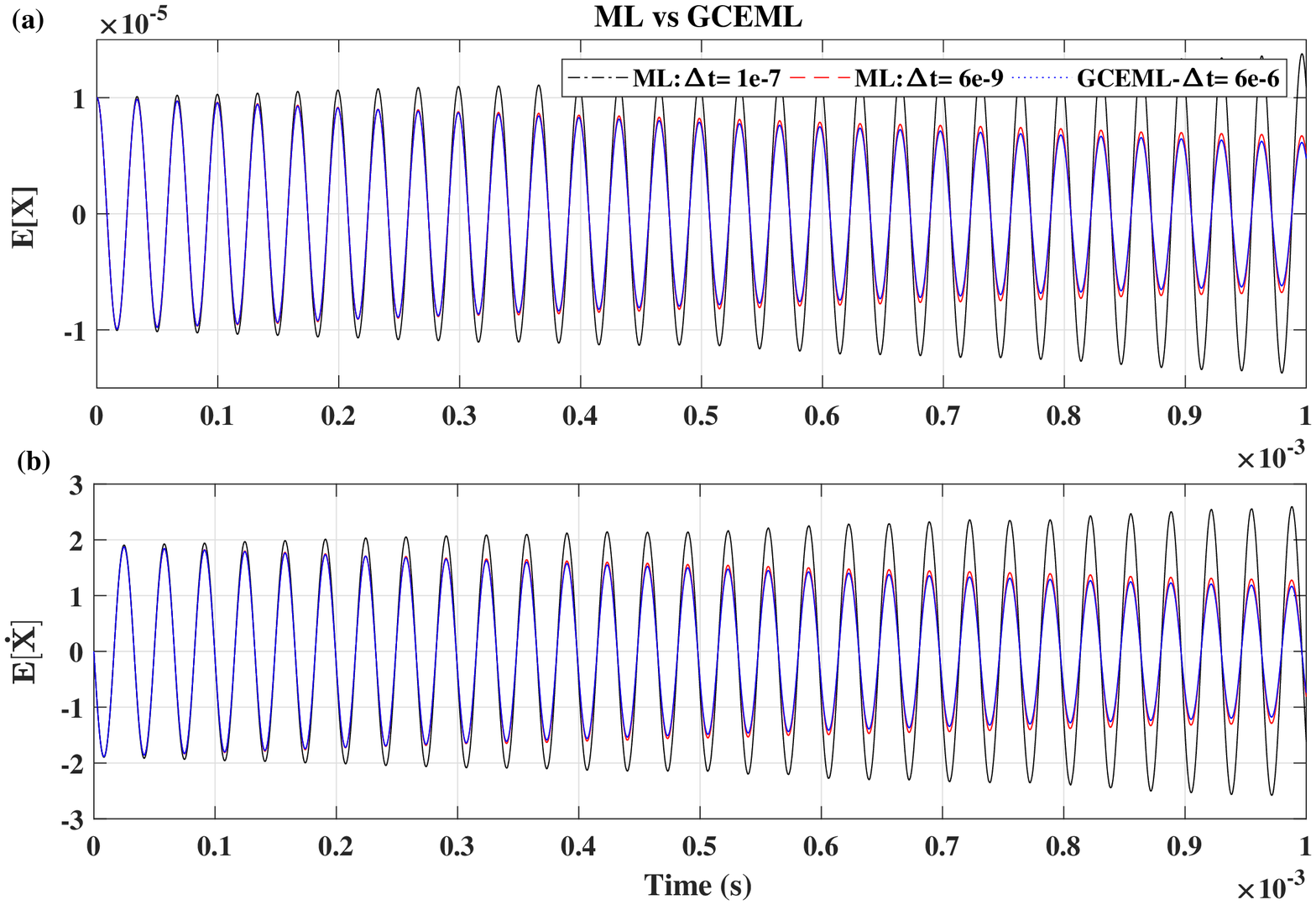}
	\caption{Ring-Type MEMS	Gyroscopes (Eq. (\ref{gyroscope})) \& GCEML: (a) ensemble mean displacement response and (b) ensemble mean velocity stability.}
	\label{figexpgyro}
\end{figure}
\begin{figure}[!htb]
	\centering
	\includegraphics[width= \textwidth]{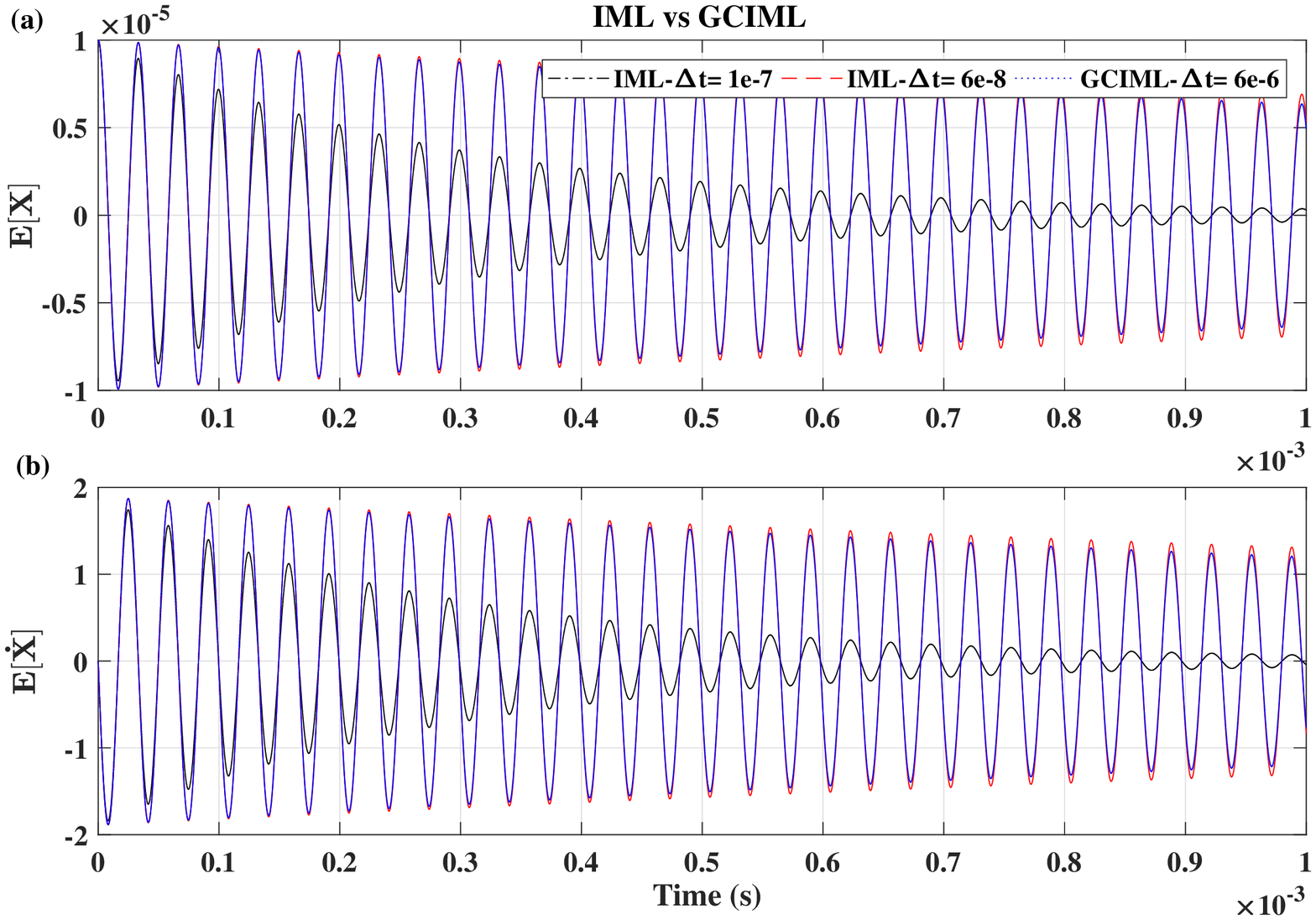}
	\caption{Ring-Type MEMS	Gyroscopes (Eq. (\ref{gyroscope})) \& GCIML: (a) ensemble mean displacement response and (b) ensemble mean velocity stability.}
	\label{figimpgyro}
\end{figure}
The properties of the ring are, $\rho$= 8800 $kg/m^3$, $E$= 210 $\times 10^9$ $N/m^2$, $r$= 500 $\mu m$, $\tilde{h}$= 12.5 $\mu m$ and $\tilde{b}$= 12.5 $\mu m$. During response simulations it is considered that, the input angular
rate $\Omega_0$ increases smoothly from $0 \to 2\pi$ rad/s between $t \in [0,0.005]s$. An initial displacement of $1 \times 10^{-5}$m is considered for the driving coordinate $q_1$. The external force is modeled using $F=6$ and $\omega=2\pi$ rad/s. The simulation is performed for $\xi$=0.8\% and $\mu_0=14.9\times 10^{-4}$. The total period of simulation is considered as 0.001s in this study. 

\begin{figure}[!htb]
	\centering
	\includegraphics[width= \textwidth]{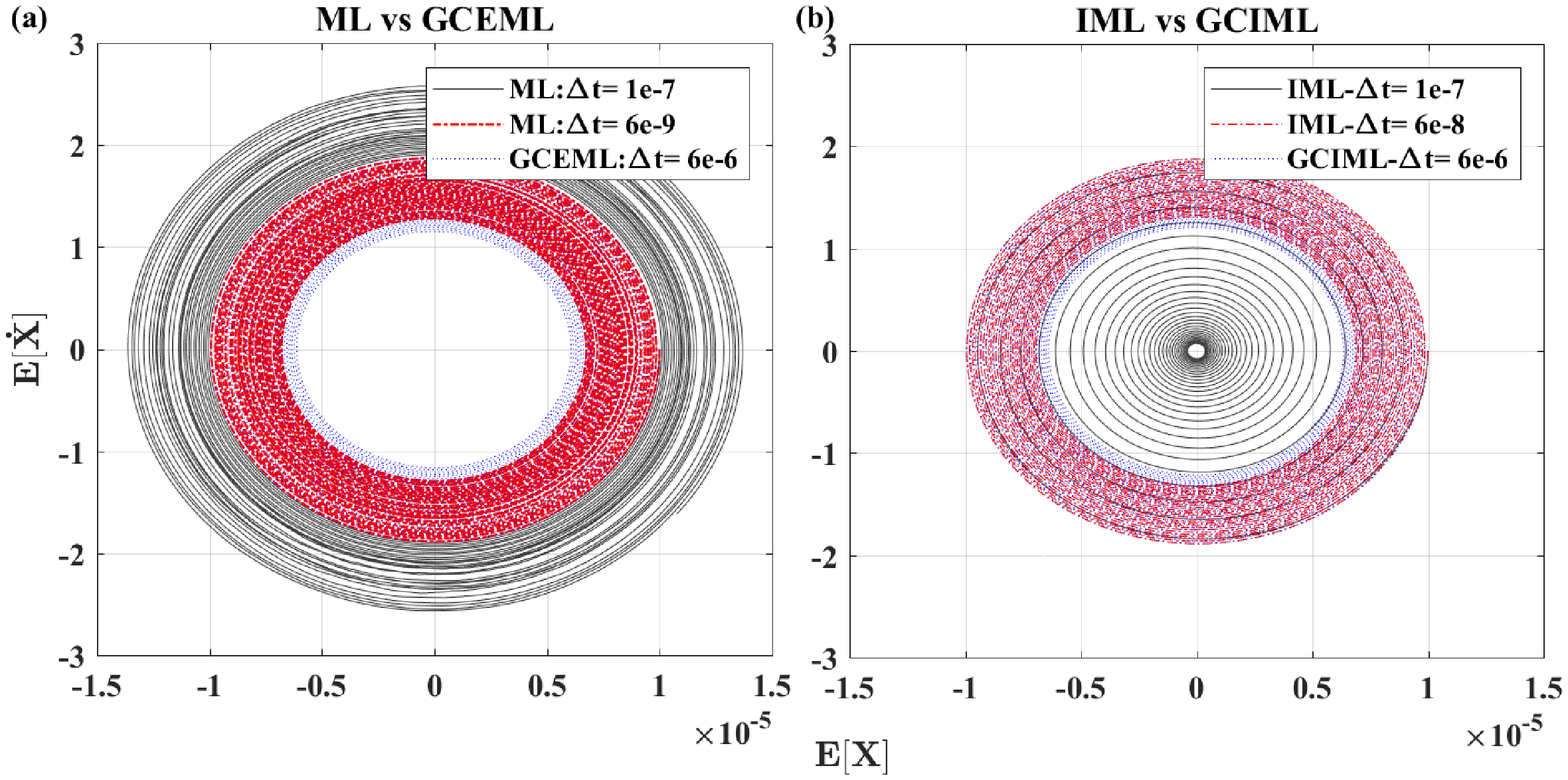}
	\caption{Ensemble mean stability - Ring-Type MEMS	Gyroscopes (Eq. (\ref{gyroscope})): (a) GCEML based solution and (b) GCIML based solution}
	\label{figppgyro}
\end{figure}

For the purposes of verifying the system state predictions, in the absence of an exact analytical solution the time responses generated by the traditional and proposed weakly corrected schemes are compared with a reference solution estimated using explicit $\Delta t$=$6 \times 10^{-9}$ and implicit $\Delta t$=$6 \times 10^{-8}$ (since the numerical convergence is claimed to reached in more than 150,000 steps in literature \cite{asokanthan2017stability}).

The time response of the MEMS oscillator using proposed GCEML and GCIML scheme against the traditional explicit and implicit schemes are displayed in Figs. \ref{figexpgyro} and \ref{figimpgyro}. For the brevity results of GCSIML is not illustrated, since the performance of GCSIML is found to approximately identical to GCIML. A number of 200 MC particles are used in measurement correction. In fig. \ref{figexpgyro}, it is evident that the ML scheme works poorly even when a slightest coarser $\Delta t$=1e-7 is adopted and predicts response which tends to diverge out from the reference prediction as the time integration progresses. In case of the corresponding alternative GCEML scheme, it is found to be robust in pathwise approximation of the reference solution with a coarser $\Delta t$=6e-6 ($>$ 1e-7). Figure \ref{figimpgyro} provides a comparison of the time responses predicted by the implicit Milstein and GCIML which reveals that, the implicit Milstein scheme tend to under-predict the response for $\Delta t$=1e-8, that introduces an increase in system decay causing divergence from the reference solution. On the contrary, the contribution of the weak correction in GCIML is significant enough to predict reference solution with $\Delta t$=6e-6. In the contrast, it is to be understood that although the order of reduction in $\Delta t \approxeq$ 5.94e-6 is small in magnitude, however the number of steps required to reach the numerical convergence of 15000 is only $\approxeq$1500 in case of the proposed GCEML, GCSIML and GCIML schemes which is very significant (reduction in 13500 steps) and of keen interest in view of study of such gyroscope oscillators.

Further, in the stability analysis it is visible in the Fig. \ref{figppgyro} that for a coarse time step the solution of GCEML and GCIML always stay in the adjacent to the reference orbit, whereas the solutions of explicit and implicit schemes without correction either form an over-predicted orbit or a spiral type center orbit, respectively. Thus, irrespective of the multi-periodic and chaotic characteristics, under the given initial and forcing condition the proposed schemes are able to correctly portray the ensemble mean stability.

\section{Conclusions}
The proposed three weakly corrected explicit, semi-implicit and implicit Milstein schemes offer computational advantage over the classical Ito-Taylor schemes, by overcoming the requirement of mathematically exhaustive formulation involving higher order MSIs in the higher-order Ito-Taylor approximations and reducing the computational demand of the lower order schemes due to the necessity of finer time sampling. The numerical case study provided the evidence towards successful implementation of the weak-correction in additive form, which is in line with the other schemes based on change of measures. From the numerical studies undertaken in this paper, it is evident that the semi-implicit and implicit Milstein approximations are more robust than the explicit one for accurate estimation of the dynamical response under stochastic excitation. The proposed approach, however, further improves the performance of the original schemes at coarser sampling rates. Thus, the proposed Girsanov corrected Milstein schemes appears to be a suitable candidate to replace the classical schemes, but further studies are required in order to understand the efficacy of such schemes for a wide class of stochastic nonlinear dynamic problems.

\textbf{Acknowledgements: } SC acknowledges the financial support received from IIT Delhi in form of seed grant.

\appendix

\section{Girsanov filtration}\label{filtering}
In order to derive the filtration equation consider a scalar valued function ${\Phi _t}(\bf{X})$ whose evolution in differential from is required. Noting that, the Radon-Nikodym derivative ${{\bf Z} _t}$ here, is the solution to the scalar SDE: $d{{\bf Z} _t} = {{\bf Z} _t}{\gamma \left( {t,{{\bf{X}}_t}} \right)}d{{\tilde {\bf W}}_t}$, and $d{\Phi _t}({\bf{X}}) = {\partial _x}\left( {{\Phi _t}({\bf{X}})} \right)d{{\bf{X}}_t} + \frac{1}{2}\partial _x^2\left( {{\Phi _t}({\bf{X}})} \right)d{\bf{X}}_t^2$, consider the following stochastic integration by parts formula \cite{calin2015informal}:
\begin{equation}\label{byparts}
\small
d\left( {{\Phi _t}\left( {\bf{X}} \right).{{\bf{Z}}_t}} \right) = {{\bf{Z}}_t}{\partial _x}{\left( {{\Phi _t}({\bf{X}})} \right)^T}d{{\bf{X}}_t} + \frac{1}{2}{{\bf{Z}}_t}\left\langle {\partial _x^2\left( {{\Phi _t}({\bf{X}})} \right){d{\bf{X}}_t^2}} \right\rangle  + {\Phi _t}\left( {\bf{X}} \right){{\bf{Z}}_t}\gamma {\left( {s,{{\bf{X}}_s}} \right)^T}d{{{\bf{\tilde W}}}_t}
\end{equation}
Where, ${\partial _x}\left( . \right)$ and $\partial _x^2\left( . \right)$ are the first and second order partial derivative with respect to x. Utilizing $d{\bf X}_t$ as the equivalent SDE form of one of the Milstein schemes (from Eq. (\ref{Q_brownian})), and noting that $\gamma \left( {t,{{\bf{X}}_t}} \right) = {\bf \rho}^{-1}({{\bf e}_{t^*}}/dt + \frac{1}{2}{\bf f}\left( {{t^*},{{\bf{X}}_{{t^*}}}} \right) )$, the Eq. (\ref{byparts}) can be rephrased as,

\begin{equation}
\begin{array}{ll}
d\left( {{\Phi _t}\left( {\bf{X}} \right).{{\bf{Z}}_t}} \right) &= {{\bf{Z}}_t}{\partial _x}{\left( {{\Phi _t}({\bf{X}})} \right)^T}\Bigl\{ {\bf{g}}\left( {{t^*},{{\bf{X}}_{{t^*}}}} \right)dt + \frac{1}{2}{\bf{f}}\left( {{t^*},{{\bf{X}}_{{t^*}}}} \right){\bf{f}}{{\left( {{t^*},{{\bf{X}}_{{t^*}}}} \right)}^T}\left( {d{\bf{W}}_t^2 - dt} \right) + \\& {\bf{f}}\left( {{t^*},{{\bf{X}}_{{t^*}}}} \right)d{{\bf{W}}_t} \Bigr\} + 
\frac{1}{2}{{\bf{Z}}_t}{{\Phi ''}_t}\left( {\bf{X}} \right) \biggl\{ {\left( {{\bf{g}}\left( {{t^*},{{\bf{X}}_{{t^*}}}} \right)} \right)^2}d{t^2} + {\left( {{\bf{f}}\left( {{t^*},{{\bf{X}}_{{t^*}}}} \right)} \right)^2}d{\bf{W}}_t^2 + \\& \frac{1}{4}{\left( {{\bf{f}}\left( {{t^*},{{\bf{X}}_{{t^*}}}} \right){\bf{f}}{{\left( {{t^*},{{\bf{X}}_{{t^*}}}} \right)}^T}} \right)^2}{\left( {d{\bf{W}}_t^2 - dt} \right)^2} + 
2{\bf{g}}\left( {{t^*},{{\bf{X}}_{{t^*}}}} \right){\bf{f}}\left( {{t^*},{{\bf{X}}_{{t^*}}}} \right)d{{\bf{W}}_t}dt + \\& {\bf{g}}\left( {{t^*},{{\bf{X}}_{{t^*}}}} \right)\left( {{\bf{f}}\left( {{t^*},{{\bf{X}}_{{t^*}}}} \right){\bf{f}}{{\left( {{t^*},{{\bf{X}}_{{t^*}}}} \right)}^T}} \right)dt\left( {d{\bf{W}}_t^2 - dt} \right) + \\&
{\bf{f}}\left( {{t^*},{{\bf{X}}_{{t^*}}}} \right)\left( {{\bf{f}}\left( {{t^*},{{\bf{X}}_{{t^*}}}} \right){\bf{f}}{{\left( {{t^*},{{\bf{X}}_{{t^*}}}} \right)}^T}} \right)d{{\bf{W}}_t}\left( {d{\bf{W}}_t^2 - dt} \right) \biggr\} + \\& {{\bf{Z}}_t}{\Phi _t}\left( {\bf{X}} \right){\left( {{{\bf{\rho }}^{ - 1}}\left( {\frac{{{{\bf{e}}_{{t^*}}}}}{{dt}} + \frac{1}{2}{\bf{f}}\left( {{t^*},{{\bf{X}}_{{t^*}}}} \right)} \right)} \right)^T}d{{{\bf{\tilde W}}}_t}
\end{array}
\end{equation}
The argument $t=t^*$ here denotes the reference point in a time discretization $\Delta t$=$(t_i - t_{i-1})$, to be considered based on the integration scheme, as per section \ref{milstein_schemes}. Further, to treat the product of the MSIs in the above equation the following identities are invoked: $d{t^2} = 0$, $d{\bf W}_t^2 = dt$, $d{{\bf W}_t}dt = 0$, $dt(d{\bf W}_t^2 - dt) = 0$, and $d{{\bf W}_t}(d{\bf W}_t^2 - dt) = 0$. The first three is due the quadratic covariation while the later two are viewed as cubic variation identities, the proof for the later two are given in \textbf{Proposition \ref{pf1}}. Then it follows,
\begin{equation}\label{diff_filter}
\small
\begin{array}{ll}
d\left( {{\Phi _t}\left( {\bf{X}} \right){{\bf{Z}}_t}} \right) &= {{\bf{Z}}_t}\biggl\{ {\partial _x}{\left( {{\Phi _t}({\bf{X}})} \right)^T}\Bigl\{ {\bf{g}}\left( {{t^*},{{\bf{X}}_{{t^*}}}} \right)dt + \frac{1}{2}{\bf{f}}\left( {{t^*},{{\bf{X}}_{{t^*}}}} \right){\bf{f}}{\left( {{t^*},{{\bf{X}}_{{t^*}}}} \right)^T}\left( {d{\bf{W}}_t^2 - dt} \right)\\&  + {\bf{f}}\left( {{t^*},{{\bf{X}}_{{t^*}}}} \right)d{{\bf{W}}_t} \Bigr\}  + \frac{1}{2}\sum\limits_{j,k = 1}^n {\sum\limits_{l = 1}^m {{{\bf{f}}^{jl}}\left( {{t^*},{{\bf{X}}_{{t^*}}}} \right){{\bf{f}}^{kl}}\left( {{t^*},{{\bf{X}}_{{t^*}}}} \right)\partial _x^2\left( {{\Phi _t}} \right)} } dt + \\& {\Phi _t}\left( {\bf{X}} \right){\left( {{{\bf{\rho }}^{ - 1}}\left( {\frac{{{{\bf{e}}_{{t^*}}}}}{{dt}} + \frac{1}{2}{\bf{f}}\left( {{t^*},{{\bf{X}}_{{t^*}}}} \right)} \right)} \right)^T}d{{{\bf{\tilde W}}}_t}\biggr\}
\end{array} 
\end{equation}
The integral form of the differential in Eq. (\ref{diff_filter}) is then formalized as,
\begin{equation}
\small
\begin{array}{ll}
{\Phi _t}\left( {\bf{X}} \right).{{\bf{Z}}_t} &= {\Phi _{{t_{i - 1}}}}\left( {\bf{X}} \right).{{\bf{Z}}_{{t_{i - 1}}}} + \int_{{t_{i - 1}}}^{{t_i}} {{\bf{Z}}_s}\biggl\{ {\partial _x}{{\left( {{\Phi _s}({\bf{X}})} \right)}^T} \Bigl\{ {\bf{g}}\left( {{s^*},{{\bf{X}}_{{s^*}}}} \right)ds + \\& \frac{1}{2}{\bf{f}}\left( {{s^*},{{\bf{X}}_{{s^*}}}} \right){\bf{f}}{{\left( {{s^*},{{\bf{X}}_{{s^*}}}} \right)}^T}\left( {d{\bf{W}}_s^2 - ds} \right) + {\bf{f}}\left( {{s^*},{{\bf{X}}_{{s^*}}}} \right)d{{\bf{W}}_s} \Bigr\}  + \\& {\Phi _s}\left( {\bf{X}} \right){{\left( {{{\bf{\rho }}^{ - 1}}\left( {\frac{{{{\bf{e}}_s}}}{{ds}} + \frac{1}{2}{\bf{f}}\left( {{s^*},{{\bf{X}}_{{s^*}}}} \right)} \right)} \right)}^T}d{{{\bf{\tilde W}}}_s} + \\& \frac{1}{2}\sum\limits_{j,k = 1}^n {\sum\limits_{l = 1}^m {{{\bf{f}}^{jl}}\left( {{s^*},{{\bf{X}}_{{s^*}}}} \right){{\bf{f}}^{kl}}\left( {{s^*},{{\bf{X}}_{{s^*}}}} \right)\partial _x^2\left( {{\Phi _s}} \right)} } ds \biggr\}
\end{array}  
\end{equation}
Taking conditional expectation with respect to the filtration $\mathcal{F}_t^m$ generated by the $\mathbb{R}^m$-valued process ${\gamma \left( {s,{{\bf{X}}_s}}\right)}$ under $Q$-measure ($\Omega$,$P$), the following can be obtained,
\begin{multline}
\small
\begin{array}{ll}
{{\rm E}_Q}\left[ {{\Phi _t}\left( {\bf X} \right){{\bf Z} _t}\left| {\mathcal{F}_t^m} \right.} \right] = {{\rm E}_Q}\left[ {{\Phi _{{t_{i - 1}}}}\left( {\bf X} \right).{{\bf Z} _{{t_{i - 1}}}}\left| {\mathcal{F}_t^m} \right.} \right] + \\ \qquad {{\rm{E}}_Q}\left[ {\int_{{t_{i - 1}}}^{{t_i}} {{{\bf{Z}}_s}{\partial _x}{{\left( {{\Phi _s}({\bf{X}})} \right)}^T}\left( {{\bf{g}}\left( {{s^*},{{\bf{X}}_{{s^*}}}} \right)ds} \right)\left| {{\cal F}_t^m} \right.} } \right] + \\ \qquad
{{\rm{E}}_Q}\left[ {\int_{{t_{i - 1}}}^{{t_i}} {\frac{1}{2}{{\bf{Z}}_s}{\partial _x}{{\left( {{\Phi _s}({\bf{X}})} \right)}^T}{\bf{f}}\left( {{s^*},{{\bf{X}}_{{s^*}}}} \right){\bf{f}}{{\left( {{s^*},{{\bf{X}}_{{s^*}}}} \right)}^T}\left( {d{\bf{W}}_s^2 - ds} \right)\left| {{\cal F}_t^m} \right.} } \right] + \\ \qquad {{\rm{E}}_Q}\left[ {\int_{{t_{i - 1}}}^{{t_i}} {{{\bf{Z}}_s}{\partial _x}{{\left( {{\Phi _s}({\bf{X}})} \right)}^T}{\bf{f}}\left( {{s^*},{{\bf{X}}_{{s^*}}}} \right)d{{\bf{W}}_s}\left| {{\cal F}_t^m} \right.} } \right] +\\ \qquad
{{\rm{E}}_Q}\left[ {\int_{{t_{i - 1}}}^{{t_i}} {{{\bf{Z}}_s}\left( {\frac{1}{2}\sum\limits_{j,k = 1}^n {\sum\limits_{l = 1}^m {{{\bf{f}}^{jl}}\left( {{s^*},{{\bf{X}}_{{s^*}}}} \right){{\bf{f}}^{kl}}\left( {{s^*},{{\bf{X}}_{{s^*}}}} \right)\partial _x^2\left( {{\Phi _s}} \right)} } } \right)ds\left| {{\cal F}_t^m} \right.} } \right] + \\ \qquad {{\rm{E}}_Q}\left[ {\int_{{t_{i - 1}}}^{{t_i}} {{{\bf{Z}}_s}{\Phi _s}\left( {\bf{X}} \right){{\left( {{{\bf{\rho }}^{ - 1}}\left( {\frac{{{{\bf{e}}_s}}}{{ds}} + \frac{1}{2}{\bf{f}}\left( {{s^*},{{\bf{X}}_{{s^*}}}} \right)} \right)} \right)}^T}d{{{\bf{\tilde W}}}_s}\left| {{\cal F}_t^m} \right.} } \right]
\end{array}
\end{multline}
Applying Fubini's theorem the expectation and integral operators can be interchanged to define ${\sigma _t}\left( . \right) = {{\rm E}_Q}\left[ {\left( . \right){{\bf Z} _t}\left| {\mathcal{F}_t^m} \right.} \right]$. Further noting the results,
\begin{itemize}
	\item $\int_{{t_{i - 1}}}^{{t_i}} {{\sigma _t}\left( {\frac{1}{2}{\partial _x}{{\left( {{\Phi _s}({\bf{X}})} \right)}^T}{\bf{f}}\left( {{s^*},{{\bf{X}}_{{s^*}}}} \right){\bf{f}}{{\left( {{s^*},{{\bf{X}}_{{s^*}}}} \right)}^T}d{\bf{W}}_s^2} \right)}  = 0$ 
	\item $\int_{{t_{i - 1}}}^{{t_i}} {{\sigma _t}\left( {{\partial _x}{{\left( {{\Phi _s}({\bf{X}})} \right)}^T}{\bf{f}}\left( {{s^*},{{\bf{X}}_{{s^*}}}} \right)d{{\bf{W}}_s}} \right)}  = 0$
\end{itemize}
the following expression can be obtained,
\begin{equation}\label{coditional}
\small
\begin{array}{ll}
{\sigma _t}\left( {\Phi ({\bf{X}})} \right) &= {\sigma _{{t_{i - 1}}}}\left( {\Phi ({\bf{X}})} \right) + \int_{{t_{i - 1}}}^{{t_i}} {\sigma _s} \Bigl\{ {\partial _x}{{\left( {{\Phi _s}({\bf{X}})} \right)}^T}{\bf{g}}\left( {{s^*},{{\bf{X}}_{{s^*}}}} \right) + \\& \frac{1}{2}\sum\limits_{j,k = 1}^n {\sum\limits_{l = 1}^m {{{\bf{f}}^{jl}}\left( {{s^*},{{\bf{X}}_{{s^*}}}} \right){{\bf{f}}^{kl}}\left( {{s^*},{{\bf{X}}_{{s^*}}}} \right)\partial _x^2\left( {{\Phi _s}} \right)} } \Bigr\} ds  + \\& \int_{{t_{i - 1}}}^{{t_i}} {{\sigma _s}\left( {{\partial _x}{{\left( {{\Phi _s}({\bf{X}})} \right)}^T}\frac{1}{2}{\bf{f}}\left( {{s^*},{{\bf{X}}_{{s^*}}}} \right){\bf{f}}{{\left( {{s^*},{{\bf{X}}_{{s^*}}}} \right)}^T}} \right)ds}  + \\& \int_{{t_{i - 1}}}^{{t_i}} {{\sigma _s}\left( {{\Phi _s}\left( {\bf{X}} \right){{\left( {{{\bf{\rho }}^{ - 1}}\left( {\frac{{{{\bf{e}}_s}}}{{ds}} + \frac{1}{2}{\bf{f}}\left( {{s^*},{{\bf{X}}_{{s^*}}}} \right)} \right)} \right)}^T}} \right)d{{{\bf{\tilde W}}}_s}} 
\end{array}
\end{equation}
For the ease of representation, three operators are introduced,
\begin{equation}
\small
\begin{array}{l}
\Im _t^0\left( {\Phi ({\bf{X}})} \right) = {\partial _x}{\left( {{\Phi _t}({\bf{X}})} \right)^T}{\bf{g}}\left( {{t^*},{{\bf{X}}_{{t^*}}}} \right) + \frac{1}{2}\sum\limits_{j,k = 1}^n {\sum\limits_{l = 1}^m {{{\bf{f}}^{jl}}\left( {{t^*},{{\bf{X}}_{{t^*}}}} \right){{\bf{f}}^{kl}}\left( {{t^*},{{\bf{X}}_{{t^*}}}} \right)} \partial _x^2\left( {\Phi _t} \right)} \\
\Im _t^1\left( {\Phi ({\bf{X}})} \right) = \frac{1}{2}{\partial _x}{\left( {{\Phi _t}({\bf{X}})} \right)^T}{\bf{f}}\left( {{t^*},{{\bf{X}}_{{t^*}}}} \right){\bf{f}}{\left( {{t^*},{{\bf{X}}_{{t^*}}}} \right)^T}\\
\Im _t^2\left( {\Phi ({\bf{X}})} \right) = {\Phi _t}\left( {\bf{X}} \right){\left( {{{\bf{\rho }}^{ - 1}}\left( {\frac{{{{\bf{e}}_{{t^*}}}}}{{dt}} + \frac{1}{2}{\bf{f}}\left( {{t^*},{{\bf{X}}_{{t^*}}}} \right)} \right)} \right)^T}
\end{array}
\end{equation} 
The first one is also known as backward Kolmogorov operator. Utilizing the above operators, the Eq. (\ref{coditional}) can be rephrased as,
\begin{equation}
\small
\begin{array}{ll}
{\sigma _t}\left( {\Phi ({\bf{X}})} \right) =& {\sigma _{{t_{i - 1}}}}\left( {\Phi ({\bf{X}})} \right) + \int_{{t_{i - 1}}}^{{t_i}} {{\sigma _s}\left( {\Im _t^0\left( {\Phi ({\bf{X}})} \right)} \right)ds}  + \int_{{t_{i - 1}}}^{{t_i}} {{\sigma _s}\left( {\Im _t^1\left( {\Phi ({\bf{X}})} \right)} \right)ds}  + \\& \int_{{t_{i - 1}}}^{{t_i}} {{\sigma _s}\left( {\Im _t^2\left( {\Phi ({\bf{X}})} \right)} \right)d{{{\bf{\tilde W}}}_s}} 
\end{array}
\end{equation}
The equivalent differential form is given as: 
\begin{equation}
d{\sigma _t}\left( \Phi ({\bf{X}}) \right) =\left( {\sigma _t}\Im _t^0\left( \Phi ({\bf{X}}) \right)+ {\sigma _t}\Im _t^1\left( \Phi ({\bf{X}}) \right) \right)dt + {\sigma _t}\Im _t^2\left( \Phi ({\bf{X}}) \right)d{\tilde {\bf{W}}_t}
\end{equation}
To invoke the change of measure $P \to Q$, a normalized function ${\pi _t}({\Phi _t}({\bf{X}}))$ which normalizes the scalar valued function ${\Phi _t}({\bf{X}})$ is considered, and denoted by ${\pi _t}\left( \Phi ({\bf{X}}) \right) = {{{\sigma _t}\left( \Phi ({\bf{X}}) \right)} \mathord{\left/
		{\vphantom {{{\sigma _t}\left( \Phi ({\bf{X}}) \right)} {{\sigma _t}\left( 1 \right)}}} \right.
		\kern-\nulldelimiterspace} {{\sigma _t}\left( 1 \right)}}$, where $	{\sigma _t}\left( \Phi ({\bf{X}}) \right) = {E_Q}\left[ {\left. {{\Phi _{\bf X}}{{\bf Z} _t}} \right|\mathcal{F}_t^m} \right]$ and ${\sigma _t}\left( 1 \right) = {E_Q}\left[ {\left. {{{\bf Z} _t}} \right|\mathcal{F}_t^m} \right]$. Then from the product rule: 
\begin{equation}
\small
d{\pi _t}\left( \Phi ({\bf{X}}) \right) = \cfrac{1}{{{\sigma _t}\left( 1 \right)}}.d{\sigma _t}\left( \Phi ({\bf{X}}) \right) + {\sigma _t}\left( \Phi ({\bf{X}}) \right)d\left( {\cfrac{1}{{{\sigma _t}\left( 1 \right)}}} \right) + d\left\langle {{\sigma _t}\left( \Phi ({\bf{X}}) \right),\left. {\cfrac{1}{{{\sigma _t}\left( 1 \right)}}} \right\rangle } \right.
\end{equation} 
where, the last term is the quadratic covariation term \cite{klebaner2005introduction}. For $\Phi ({\bf{X}})=1$, it can found that, $d{\sigma _t}\left( 1 \right) = {\sigma _t}{\gamma\left( {s,{{\bf{X}}_s}} \right)}^Td{\tilde {\bf{W}}_s}$. Noting, ${\pi _t}\left( \Phi  \right) = \frac{{{\sigma _t}\left( \Phi  \right)}}{{{\sigma _t}\left( 1 \right)}}$ from the Ito's lemma $d\left( {\frac{1}{{{\sigma _t}\left( 1 \right)}}} \right)$ can be found as, 
\begin{equation}
\small
\begin{array}{ll}
d\left( {\frac{1}{{{\sigma _t}\left( 1 \right)}}} \right) &=  - \frac{1}{{{\sigma _t}\left( 1 \right)}}{\pi _t}{\left( {{{\bf{\rho }}^{ - 1}}\left( {\frac{{{{\bf{e}}_{{t^*}}}}}{{dt}} + \frac{1}{2}{\bf{f}}\left( {{t^*},{{\bf{X}}_{{t^*}}}} \right)} \right)} \right)^T}d{{{\bf{\tilde W}}}_t} + \\& \frac{1}{{{\sigma _t}\left( 1 \right)}}{\pi _t}{\left( {{{\bf{\rho }}^{ - 1}}\left( {\frac{{{{\bf{e}}_{{t^*}}}}}{{dt}} + \frac{1}{2}{\bf{f}}\left( {{t^*},{{\bf{X}}_{{t^*}}}} \right)} \right)} \right)^T}{\pi _t}\left( {{{\bf{\rho }}^{ - 1}}\left( {\frac{{{{\bf{e}}_{{t^*}}}}}{{dt}} + \frac{1}{2}{\bf{f}}\left( {{t^*},{{\bf{X}}_{{t^*}}}} \right)} \right)} \right)dt
\end{array}
\end{equation}
Then applying the quadratic identities, $d{t^2} = 0$, $dW_t^2 = dt$, and $d{W_t}dt = 0$, the following result is obtained,
\begin{equation}
\small
\begin{array}{l}
d{\pi _t}\left( \Phi  \right) = {\pi _t}\Bigl( {{\partial _x}{{\left( {{\Phi _t}} \right)}^T}{\bf{g}}\left( {{t^*},{{\bf{X}}_{{t^*}}}} \right) + \frac{1}{2}\sum\limits_{j,k = 1}^n {\sum\limits_{l = 1}^m {{{\bf{f}}^{jl}}\left( {{t^*},{{\bf{X}}_{{t^*}}}} \right){{\bf{f}}^{kl}}\left( {{t^*},{{\bf{X}}_{{t^*}}}} \right)\partial _x^2\left( {{\Phi _t}} \right)} } } \Bigr)dt + \\
{\pi _t}\left( {{\Phi _t}{{\left( {{{\bf{\rho }}^{ - 1}}\left( {\frac{{{{\bf{e}}_{{t^*}}}}}{{dt}} + \frac{1}{2}{\bf{f}}\left( {{t^*},{{\bf{X}}_{{t^*}}}} \right)} \right)} \right)}^T}} \right)d{{{\bf{\tilde W}}}_t} - {\pi _t}\left( \Phi  \right){\pi _t}{\left( {{{\bf{\rho }}^{ - 1}}\left( {\frac{{{{\bf{e}}_{{t^*}}}}}{{dt}} + \frac{1}{2}{\bf{f}}\left( {{t^*},{{\bf{X}}_{{t^*}}}} \right)} \right)} \right)^T}d{{{\bf{\tilde W}}}_t}\\
+ {\pi _t}\left( \Phi  \right){\pi _t}{\left( {{{\bf{\rho }}^{ - 1}}\left( {\frac{{{{\bf{e}}_{{t^*}}}}}{{dt}} + \frac{1}{2}{\bf{f}}\left( {{t^*},{{\bf{X}}_{{t^*}}}} \right)} \right)} \right)^T}{\pi _t}\left( {{{\bf{\rho }}^{ - 1}}\left( {\frac{{{{\bf{e}}_{{t^*}}}}}{{dt}} + \frac{1}{2}{\bf{f}}\left( {{t^*},{{\bf{X}}_{{t^*}}}} \right)} \right)} \right)dt\\
- {\pi _t}\left\{ {\Phi {{\left( {{{\bf{\rho }}^{ - 1}}\left( {\frac{{{{\bf{e}}_{{t^*}}}}}{{dt}} + \frac{1}{2}{\bf{f}}\left( {{t^*},{{\bf{X}}_{{t^*}}}} \right)} \right)} \right)}^T}} \right\}{\pi _t}{\left( {{{\bf{\rho }}^{ - 1}}\left( {\frac{{{{\bf{e}}_{{t^*}}}}}{{dt}} + \frac{1}{2}{\bf{f}}\left( {{t^*},{{\bf{X}}_{{t^*}}}} \right)} \right)} \right)^T}dt\\
+ \frac{1}{2}{\pi _t}\left( {{\partial _x}{{\left( {{\Phi _t}} \right)}^T}{\bf{f}}\left( {{t^*},{{\bf{X}}_{{t^*}}}} \right){\bf{f}}{{\left( {{t^*},{{\bf{X}}_{{t^*}}}} \right)}^T}} \right)dt
\end{array}
\end{equation}
In terms of the operators, the above equation can however be written in reduced form:
\begin{equation}
\small
\begin{array}{l}
d{\pi _t}\left( \Phi  \right) = {\pi _t}\left( {\Im _t^1\left( \Phi  \right)} \right)dt + \left( {{\pi _t}\left( {\Im _t^2\left( \Phi  \right)} \right) - {\pi _t}\left( \Phi  \right){\pi _t}{{\left( {{{\bf{\rho }}^{ - 1}}\left( {\frac{{{{\bf{e}}_{{t^*}}}}}{{dt}} + \frac{1}{2}{\bf{f}}\left( {{t^*},{{\bf{X}}_{{t^*}}}} \right)} \right)} \right)}^T}} \right)d{{{\bf{\tilde W}}}_t}\\
- \left\{ {{\pi _t}\left( {\Im _t^2\left( \Phi  \right)} \right) - {\pi _t}\left( \Phi  \right){\pi _t}{{\left( {{{\bf{\rho }}^{ - 1}}\left( {\frac{{{{\bf{e}}_{{t^*}}}}}{{dt}} + \frac{1}{2}{\bf{f}}\left( {{t^*},{{\bf{X}}_{{t^*}}}} \right)} \right)} \right)}^T}} \right\}{\pi _t}\left( {{{\bf{\rho }}^{ - 1}}\left( {\frac{{{{\bf{e}}_{{t^*}}}}}{{dt}} + \frac{1}{2}{\bf{f}}\left( {{t^*},{{\bf{X}}_{{t^*}}}} \right)} \right)} \right)dt\\
+ {\pi _t}\left( {\Im _t^1\left( \Phi  \right)} \right)dt
\end{array}
\end{equation}
After suitable substitutions the final filtering equation can be derived as,
\begin{multline}
\small
\begin{array}{ll}
d{\pi _t}\left( \Phi  \right) =& {\pi _t}\left( {\Im _t^0\left( \Phi  \right) + \Im _t^1\left( \Phi  \right)} \right)dt + \\& \left( {{\pi _t}\left( {\Im _t^2\left( \Phi  \right)} \right) - {\pi _t}\left( {{\Phi _t}\left( {\bf{X}} \right)} \right){\pi _t}{{\left( {{{\bf{\rho }}^{ - 1}}\left( {\frac{{{{\bf{e}}_{{t^*}}}}}{{dt}} + \frac{1}{2}{\bf{f}}\left( {{t^*},{{\bf{X}}_{{t^*}}}} \right)} \right)} \right)}^T}} \right)d{\psi _t}
\end{array}
\end{multline}
where, $d{I_t} = \left( {d{{{\bf{\tilde W}}}_t} - {\pi _t}\left( {{{\bf{\rho }}^{ - 1}}\left( {\frac{{{{\bf{e}}_{{t^*}}}}}{{dt}} + \frac{1}{2}{\bf{f}}\left( {{t^*},{{\bf{X}}_{{t^*}}}} \right)} \right)} \right)dt} \right)$ is the innovation term. One can find the integral representation as:
\begin{equation}
\begin{array}{ll}
{\pi _t}\left( \Phi  \right) =& {\pi _{{t_{i - 1}}}}\left( \Phi  \right) + \int_{{t_{i - 1}}}^{{t_i}} {\left( {{\pi _s}\left( {\Im _s^0\left( \Phi  \right)} \right) + {\pi _s}\left( {\Im _s^1\left( \Phi  \right)} \right)} \right)ds}  + \\& \int_{{t_{i - 1}}}^{{t_i}} {\left( {{\pi _s}\left( {\Im _s^2\left( \Phi  \right)} \right) - {\pi _s}\left( {{\Phi _s}} \right){\pi _s}{{\left( {\gamma \left( {s,{{\bf{X}}_s}} \right)} \right)}^T}} \right)d{\psi _s}} 
\end{array}
\end{equation}
The choice of $\gamma \left( {s,{{\bf{X}}_s}} \right) = {{\bf{\rho }}^{ - 1}}\left( {\frac{{{{\bf{e}}_{{t^*}}}}}{{dt}} + \frac{1}{2}{\bf{f}}\left( {{t^*},{{\bf{X}}_{{t^*}}}} \right)} \right)$, shall be done purely based on the type of Milstein scheme in interest.


\end{document}